%% file: draft_v1.tex
\documentclass[preprint,12pt]{elsarticle}

\date{}

\usepackage{caption}
\usepackage{subcaption}
\usepackage{multirow}
\usepackage{amsmath,amssymb,bm}
\usepackage{algorithm}
\usepackage{algpseudocode} 
\usepackage{xcolor}
\usepackage{float}
\usepackage{placeins} 
\usepackage[a4paper,margin=2.5cm]{geometry}
\usepackage{tikz}
\usetikzlibrary{positioning, arrows.meta}
\usepackage{booktabs}
\usepackage{graphicx}
\usepackage{caption}
\usepackage{tabularx} 
\usepackage{lmodern}
\usepackage{fix-cm}
\usetikzlibrary{positioning, shapes, arrows.meta, fit, backgrounds, calc}

\usepackage{amssymb}
\usepackage{amsmath}
\usepackage{tikz}
\usepackage{rotating}
\usepackage{hyperref}
\usepackage{amsmath, amssymb}

\usepackage{tikz}
\usetikzlibrary{arrows.meta, positioning, shapes, decorations.pathreplacing, calc}
\usepackage{comment}

\usepackage{booktabs}
\usepackage{siunitx}
\usepackage[percent]{overpic}
\usepackage[utf8]{inputenc}
\usepackage{siunitx}
\usepackage{graphicx}

\usepackage{subcaption}

\definecolor{myGreen}{rgb}{0.4501, 0.6784,0.2745}

\definecolor{myblue}{rgb}{0.705882, 0.78039, 0.90588}

\definecolor{mywhite}{rgb}{1,1,1}
\definecolor{blue2}{rgb}{0.7098,0.78039,0.913725}

\definecolor{blueshape}{rgb}{0.670589,0.74902,0.901961}

\definecolor{greenframe}{rgb}{0.77255,1,0.77255}

\definecolor{greenabstract}{rgb}{0.55691,0.851,0.451}

\definecolor{blueabstract}{rgb}{0.27451,0.6941,0.88235}

\journal{Ocean Engineering}

\begin{document}

\begin{frontmatter}





\title{Localized Generative Design of Ship Hulls with Diffusion Models and High-Fidelity Assessment}

\author[label1]{Hadi Keramati\corref{cor1}}
\ead{hadi.keramati@ubc.ca}

\author[label1]{Mohammed A. Hannan}

\author[label1]{Patrick Kirchen}

\author[label2]{Sina Tousi}

\author[label2]{Joonas Riekkinen}

\author[label3]{Jim Pyra}

\author[label1]{Rajeev K. Jaiman}

\cortext[cor1]{Corresponding author: hadi.keramati@ubc.ca}

\affiliation[label1]{
    organization={Department of Mechanical Engineering, University of British Columbia},
    city={Vancouver},
    state={BC},
    country={Canada}
}

\affiliation[label2]{
    organization={Elomatic Consulting Inc.},
    city={Vancouver},
    state={BC},
    country={Canada}
}

\affiliation[label3]{
    organization={InnovMarine Inc.},
    city={L\'evis},
    state={QC},
    country={Canada}
}
\begin{abstract}
The ship design process is highly iterative, requiring repeated cycles of geometric modification, performance evaluation, and expert judgment. Recent advances in generative diffusion models show considerable promise for automating hull-form generation. However, most existing approaches focus on global, conceptual redesign and provide limited support for practical engineering workflows, which require localized refinement of specific hull features and high-fidelity performance evaluations to establish confidence in the generated designs.
This paper presents a reward-directed diffusion framework for localized ship hull optimization, enabling targeted modification of specific geometric regions while preserving global hull form continuity. The proposed framework introduces a parametric masking mechanism in dataset generation that constrains the generative process to user-selected subsets of design variables, such as the bulbous bow and forebody, allowing controlled and physically smooth geometric refinement within a full-hull representation. A surrogate-assisted sequential Monte Carlo sampling strategy guides the diffusion process toward hull geometries with reduced predicted calm-water resistance, with stability retained through the training data and displacement matched in a subsequent screening.
An automated viscous free-surface CFD pipeline, validated against model-scale measurements for two benchmark container ships, is then applied to a generated variant at the displacement of the KRISO Container Ship (KCS). The total resistance of this variant is comparable to that of KCS at $F_n = 0.195$ but about 45 to 50 percent higher at $F_n = 0.250$ and $0.260$, an increase that the empirical screening did not capture. To examine the parametric representation within which the variants are generated, the pipeline is also applied to GenHull~E, a hull refined by a naval architect toward KCS within the same representation. The refined hull has six to eight percent lower total resistance than the benchmark up to $F_n = 0.240$ and higher total resistance above $F_n = 0.250$, with the crossover occurring where a friction advantage of $8.07$\,N and a pressure penalty of $8.49$\,N nearly cancel. The two components are measured separately and have opposite speed dependence, so that an objective aggregated across the crossover can obscure the opposing performance trends.
The evaluation also shows that the representation cannot express a faired bulbous bow, and it indicates limitations of the reward formulation, which contains neither viscous pressure nor a stability term. These findings demonstrate how high-fidelity evaluation can guide improvements to design representations and objectives, helping ensure that predicted gains are physically realizable.

\end{abstract}

\end{frontmatter}


\section{Introduction}
The optimization of ship design involves a multidisciplinary and tightly coupled workflow integrating geometric modeling, hydrodynamic analysis, and numerical optimization within broader naval-architectural requirements \cite{gaspar_2012}. For hull-form development, a central objective is to reduce resistance while retaining the displacement, carrying capacity, and operational characteristics required of the vessel. Lower resistance reduces the effective power required for propulsion, but resistance reduction must be considered alongside hydrostatic, geometric, and operational requirements rather than pursued as an isolated objective. The resulting optimization problem therefore combines nonlinear performance objectives with constraints on vessel particulars, hydrostatic properties, geometric feasibility, and allowable design variations \cite{zhang2019research}. In practice, this search frequently begins from an established parent hull rather than an unrestricted design space \cite{perez2011constrained}, with selected features modified while the remainder of the vessel is retained. Han et al.~\cite{han2012hydrodynamic}, for example, optimized the forebody while preserving the aftbody geometry. Such localized optimization requires a parameterization sufficiently flexible to represent useful modifications while controlling their extent and maintaining smooth transitions to retained surfaces. Expert judgment therefore remains important in defining the design envelope and interpreting trade-offs, while automation supports systematic exploration within this engineering-defined space.

The bulbous bow and adjacent forebody are particularly relevant to localized resistance optimization because their geometry influences the bow-wave system and immersed-volume distribution. Their hydrodynamic benefit is strongly dependent on operating condition: a modification advantageous at one speed or draft may provide less benefit, or become detrimental, at another. Lu et al.~\cite{lu2016hydrodynamic} demonstrated this dependence through bulbous-bow optimization under multiple operating conditions. Candidate designs should therefore be assessed across the relevant operating envelope rather than at a single design point, while accounting for accompanying changes in displacement, draft, and other characteristics relative to the parent hull.
Evaluating these alternatives requires balancing prediction fidelity against computational cost. Semi-empirical methods, such as the Holtrop-Mennen resistance and power prediction method, provide rapid estimates based on relationships calibrated from existing ship data \cite{holtrop1982approximate}, but may not resolve the hydrodynamic effects of localized geometric changes. Viscous free-surface CFD provides a higher-fidelity assessment of geometry-dependent wave, pressure, and viscous effects, but its repeated use in many candidates and operating conditions is computationally expensive. This motivates workflows that use inexpensive methods for design-space exploration while reserving high-fidelity CFD for selected candidates \cite{zhang2019research}.

Surrogate-assisted optimization provides an established means of reducing this expense. By approximating the relationship between design variables, operating conditions, and hydrodynamic responses, a surrogate can support the screening of many candidates at a fraction of the cost of repeated CFD evaluations. Copped\`e et al.~\cite{coppede2019hydrodynamic} combined high-fidelity CFD with Gaussian-process response surfaces for KCS hull optimization and reassessed selected designs using refined CFD calculations. 
This separation between inexpensive exploration and subsequent numerical verification is important because a favorable surrogate prediction does not, by itself, establish a hydrodynamic improvement. The quality and coverage of the training data, the representation of relevant operating conditions, and the distance of proposed designs from the evaluated region all influence the reliability of the search. Consequently, surrogate-based design should be understood as a means of identifying promising candidates for further assessment, not as a substitute for the verification of the final hull forms. Even with inexpensive evaluations, the optimization efficiency depends on the geometric representation and search strategy. Gradient-based and derivative-free hull optimization depends on response smoothness, sensitivities, dimensionality, and evaluation budget. Reducing the design space to physically meaningful shape variations can therefore be as important as selecting the optimization algorithm. Harries and Abt~\cite{harries2019faster} examined shape variation and reduced design spaces to accelerate simulation-driven design. These considerations motivate the use of existing hull families to guide coordinated parameter changes beyond exploration within individual variable bounds.

Data-driven generative models provide a complementary approach to constructing design spaces by learning relationships among geometric variables and using them to generate new candidate forms. Khan et al.~\cite{khan2023shiphullgan} developed ShipHullGAN as a learned parametric representation for ship-hull generation and design-space exploration. Such models can retain correlations among hull characteristics while exploring beyond individual training examples, although geometric plausibility, hydrodynamic performance, and vessel-specific requirements must still be addressed explicitly. Diffusion models generate new samples by learning to reverse a progressive noise-corruption process \cite{sohl2015deep,ho2020denoising}. For hull design, ShipGen demonstrated diffusion-based generation with feasibility and performance guidance \cite{bagazinski2023shipgen}, while C-ShipGen introduced conditional generation using prescribed principal dimensions and guidance from a total-resistance regression model \cite{bagazinski2024c}. For design optimization, however, learning the distribution of feasible geometries is only a starting point. The sampling process must also be directed toward candidates that improve a prescribed objective while satisfying relevant constraints. Krishnamoorthy et al.~\cite{krishnamoorthy2023diffusion} demonstrated objective-guided conditional diffusion for offline black-box optimization. In hull design, such guidance can connect generative sampling to hydrodynamic performance, for example by rewarding lower predicted resistance while constraining unacceptable designs.

Performance information can be incorporated into diffusion-based design through several mechanisms. Differentiable surrogates enable gradient-based guidance, as demonstrated for resistance-directed hull generation in C-ShipGen \cite{bagazinski2024c} and coupled thermal-fluid heat-sink design in HeatGen \cite{keramati2025heatgen}. Alternatively, reward-directed methods use performance values without requiring reward derivatives. Keramati et al.~\cite{keramati2026reward} combined reward-directed fine-tuning with value-based sampling for engineering design, including ship-hull applications. These approaches differ primarily in their requirements for differentiable feedback, model retraining, and sampling-time evaluation. Inference-time guidance is particularly useful when a trained generative model is reused for different performance objectives, since sampling can be redirected without retraining the diffusion model. Li et al.~\cite{li2024derivative} demonstrated derivative-free value-based guidance without differentiating the reward or fine-tuning the generative model. Sequential Monte Carlo (SMC) provides a population-based alternative that propagates and resamples weighted particles \cite{del2006sequential}, allowing reward-guided diffusion to preferentially sample promising candidate trajectories while retaining the learned denoising model as the proposal mechanism.

In this study, we propose a diffusion framework for localized optimization of the KRISO Container Ship (KCS), focusing on the bulbous bow and adjacent forebody. A diffusion model is trained on KCS-based geometries with controlled variations in this region, while parametric masking retains the remaining hull parameters at their baseline values. A surrogate provides rapid performance evaluation, and reward-guided sequential Monte Carlo (SMC) sampling directs generation toward lower predicted resistance. Stability is retained through the training data and displacement is matched during subsequent screening, enabling data-driven exploration within an engineering-defined design envelope. The contribution is twofold. First, we introduce a derivative-free, reward-directed diffusion framework trained on a masked dataset for localized ship-hull optimization. Second, we provide a viscous hydrodynamic assessment of a sampled variant and of a hull refined by a naval architect within the same parametric representation. An automated CFD pipeline, with mesh and solver parameters determined systematically from hull geometry and Froude number, is validated against two benchmark container ships and applied to both hulls. The assessment reveals limitations in geometric representation, reward fidelity, operating-range aggregation, and preservation of stability margin that are not apparent from surrogate-based optimization alone.

The remainder of this paper is structured as follows. Section~\ref{sec:methodology} presents the diffusion framework, including the masked design space, the denoising diffusion model, the reward-guided SMC sampling strategy, and the automated CFD pipeline together with its validation against the KCS and Duisburg Test Case benchmarks. Section~\ref{sec:results} presents diagnostics of the diffusion model and sampler, summarizes the empirical resistance screening of the generated hull forms and the CFD assessment of one screened variant, and then reports the viscous assessment of the refined hull across seven Froude numbers, covering geometry, hydrostatics, resistance components, wave field and hull pressure. Section~\ref{sec:discussion} discusses the speed dependence of the resistance comparisons and the representational, evaluative and constraint limitations they expose, and Section~\ref{sec:conc} concludes.

\section{Methodology}
\label{sec:methodology}
This section presents the parametric diffusion framework for localized ship-hull optimization. The methodology enables targeted refinement of selected hull regions while preserving global hull-form continuity. It combines a denoising diffusion model with parametric masking and surrogate-guided sampling, with stability enforced through the training data and displacement matched during subsequent screening.

Fourteen design variables controlling the bow and bulbous-bow geometry are selected as the active subset, while the remaining variables are fixed at their baseline KCS values. Latin Hypercube Sampling (LHS) \cite{shields2016generalization} is used to sample the active design space, producing 3,800 designs that meet the prescribed stability constraints. The full design vectors are used to train a hydrodynamic-resistance surrogate, with labels computed following Keramati et al.~\cite{keramati2026reward}: wave-making resistance is estimated using Michell's integral and frictional resistance using the ITTC-57 correlation line applied to the wetted surface. The surrogate is a residual neural network (ResNet) with four residual blocks and 256 hidden units. Sampling is controlled through KL regularization and a temperature parameter.

\begin{figure}[h]
    \centering
    \includegraphics[width=\textwidth]{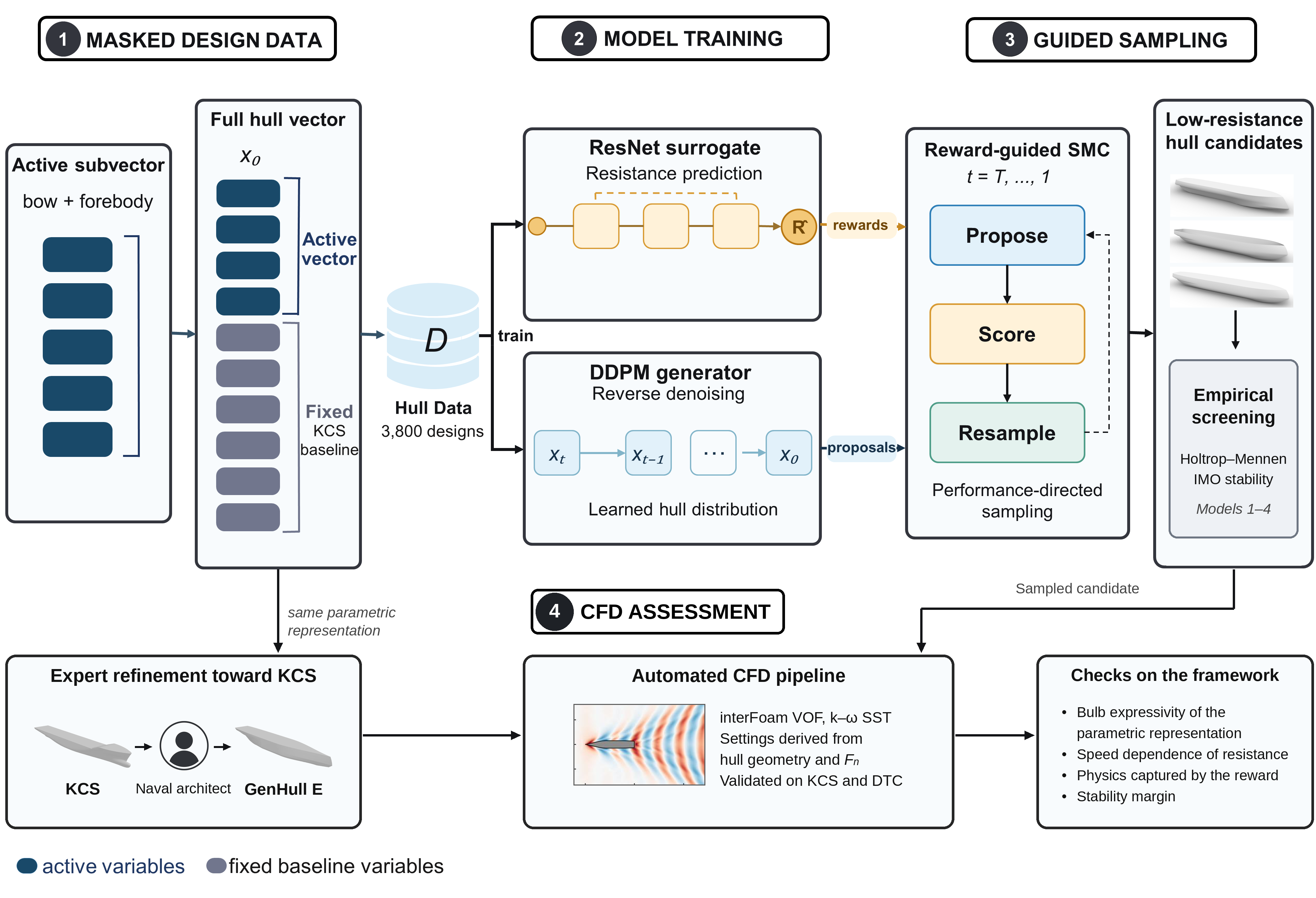} 
\caption{Schematic of the ship-hull optimization workflow with a partially variable design vector. The framework comprises data generation, training of a denoising diffusion probabilistic model (DDPM) and a deep residual surrogate model, reward-guided sampling via sequential Monte Carlo (SMC), and CFD assessment of a sampled variant and of a hull refined by a naval architect.} 
    \label{fig:schem2026} 
\end{figure}

The DDPM is trained on the fused full design vectors, while masking is imposed at the dataset level by varying only the active subset. This preserves a full-hull representation while concentrating geometric variability in the bow and bulbous-bow regions. After training, surrogate-guided sequential Monte Carlo sampling directs generation toward designs with lower predicted resistance. Figure~\ref{fig:schem2026} summarizes the workflow. In Stage~4, the CFD pipeline of Section~\ref{subsec:cfd_pipeline} is applied to Model~1, one of the sampled variants (Section~\ref{subsec:empirical_screening}), and to GenHull~E, a hull refined by a naval architect toward KCS within the same parametric representation (Section~\ref{subsec:cfd_genhull}). The DDPM and reward-guided sampling procedure are described next.

\subsection{DDPM Formulation}
\label{subsec:ddpm}

Diffusion models are likelihood-based generative models that transform data into noise via a fixed forward Markov chain and learn the reverse denoising dynamics with a neural network \cite{sohl2015deep, ho2020denoising}. Given a data sample $x_0\in\mathbb{R}^d$, DDPMs introduce latent variables $x_{1:T}=\{x_1,\dots,x_T\}$ over $T$ diffusion steps and define the generative model
\begin{equation}
p_\theta(x_{0:T}) = p(x_T)\prod_{t=1}^{T} p_\theta(x_{t-1}\mid x_t),
\qquad
p(x_T)=\mathcal{N}(0,I),
\label{eq:latent}
\end{equation}
where $\theta$ denotes the parameters of the neural network and $I$ is an identity matrix.

\paragraph{Forward (noising) process.}
A fixed variance schedule $\{\beta_t\}_{t=1}^T$, with $\beta_t\in(0,1)$, defines the forward Markov transitions
\begin{equation}
q(x_t \mid x_{t-1}) = \mathcal{N}\!\left(x_t;\sqrt{1-\beta_t}\,x_{t-1},\beta_t I\right),
\label{eq:forward}
\end{equation}
where $q(\cdot)$ denotes the forward diffusion distribution.

Defining $\alpha_t=1-\beta_t$ and $\bar{\alpha}_t=\prod_{s=1}^t \alpha_s$, the marginal distribution of $x_t$ given $x_0$ admits the closed form
\begin{equation}
x_t=\sqrt{\bar{\alpha}_t}\,x_0+\sqrt{1-\bar{\alpha}_t}\,\epsilon,
\quad \epsilon\sim\mathcal{N}(0,I),
\label{eq:marginal}
\end{equation}
where $\epsilon\in\mathbb{R}^d$ is standard Gaussian noise.

\paragraph{Reverse (denoising) process.}
The reverse-time transitions are modeled as Gaussian conditionals
\begin{equation}
p_\theta(x_{t-1}\mid x_t)=\mathcal{N}\!\left(x_{t-1};\mu_\theta(x_t,t),\sigma_t^2 I\right),
\label{eq:reverse}
\end{equation}
where $\mu_\theta(x_t,t)$ is a learned mean function and $\sigma_t^2$ denotes the variance at step $t$.

Using the noise-parameterization of \cite{ho2020denoising}, the reverse update can be written as
\begin{equation}
x_{t-1}
=
\frac{1}{\sqrt{\alpha_t}}
\Bigl(
x_t-\frac{\beta_t}{\sqrt{1-\bar{\alpha}_t}}\,\epsilon_\theta(x_t,t)
\Bigr)
+\sigma_t z,
\quad z\sim\mathcal{N}(0,I),
\label{eq:denoise}
\end{equation}
where $\epsilon_\theta(x_t,t)$ is a neural network predicting the injected noise and $z\in\mathbb{R}^d$ is standard Gaussian noise.

\paragraph{Training objective.}
Optimizing the evidence lower bound yields the standard DDPM loss
\begin{equation}
\begin{aligned}
\mathcal{L}_{\text{simple}}
&=
\mathbb{E}_{x_0,\,\epsilon,\,t}
\!\left[
\|\epsilon-\epsilon_\theta(x_t,t)\|^2
\right], \\
x_t
&=
\sqrt{\bar{\alpha}_t}\,x_0+\sqrt{1-\bar{\alpha}_t}\,\epsilon,
\quad
t\sim \mathrm{Unif}\{1,\dots,T\},
\end{aligned}
\label{eq:ddpm_final_loss}
\end{equation}
where the expectation is taken over the data distribution, Gaussian noise $\epsilon$, and uniformly sampled diffusion steps $t$.
The diffusion formulation described above provides a generative mechanism for producing physically smooth and coherent hull-form variations within the prescribed parametric design space. However, in its standard form, the reverse diffusion process samples designs according to the learned data distribution without preferential emphasis on hydrodynamic performance. To transform the diffusion model from a passive generator into an active design optimization tool, a performance-aware guidance mechanism is required. 
In the present work, this guidance is introduced through a reward-based Sequential Monte Carlo (SMC) sampling strategy, which biases the reverse diffusion process toward hull geometries with lower predicted resistance while preserving diversity and geometric continuity. The details of this reward-guided sampling procedure are described next.

\subsection{Sequential Monte Carlo sampling}
\label{subsec:smc}
 After training, the surrogate model is used to evaluate candidate designs and assign rewards corresponding to their predicted hydrodynamic performance, with higher reward indicating lower total resistance. Sequential Monte Carlo sampling, following Algorithm~\ref{alg:smc_guidance}, is then employed to bias the reverse diffusion process toward low-resistance hull geometries while remaining close to the learned data distribution. By balancing performance-driven selection with probabilistic resampling, the SMC procedure enables the diffusion model to function as an active optimization tool, promoting high-performing designs while preserving diversity and smooth geometric transitions.

In Algorithm~\ref{alg:smc_guidance}, the function $\hat v_t(x_t)$ denotes the reward surrogate that evaluates the quality of a design state at diffusion step $t$. At each step, $\hat v_t(\cdot)$ is obtained by applying the ResNet surrogate model to the denoised design representation corresponding to $x_t$ and returns a scalar reward proportional to the predicted drag value.
In this work, the reward is defined such that higher values of $\hat v_t$ correspond to lower predicted drag.
This reward applies the importance weights through an exponential tilting term $\exp(\hat v_t/\alpha)$, which biases the sampling process toward high-performing designs while preserving diversity through the temperature parameter $\alpha$. The resampling threshold is set to $0.5$ to mitigate weight degeneracy and prevent particle collapse, ensuring that the sampling process is not dominated by a single trajectory during the reverse diffusion steps.
 From an engineering perspective, this procedure biases the generative process toward lower-resistance hull forms while retaining sufficient diversity to avoid premature convergence to a single design. For many generative settings, Algorithm~\ref{alg:smc_guidance} remains unchanged, with the hyperparameter $\alpha$ treated as a domain-dependent variable. Before presenting the results and discussion in the next section, we summarize the practical steps required for industrial practitioners and engineers to deploy generative models in practice. The first step is data curation in the structured format shown in Figure \ref{fig:schem2026}. It is generally advisable to begin with the full design vector, run and monitor the system performance, and validate feasibility. If successful, a component-wise model based on subvectors can then be constructed for localized optimization. Once sufficient data are generated and the empirical distribution is not severely skewed, the dataset can be used to train both predictive and generative models. These steps concern only the generation of data consistent with the model input representation. Subsequently, Algorithm~\ref{alg:smc_guidance} can be applied with a moderate temperature value to obtain improved designs without excessively biasing the search toward extreme high-performance regions. After robustness is confirmed, a systematic search over the temperature parameter can be conducted to identify the value that yields the best-performing designs.
 
 We used Weights and Biases platform \cite{wandb} to automate the search for parameters and hyperparameters. The final values of the parameters and hyperparameters of both the DDPM and ResNet models , including the corresponding loss functions, are reported in Table \ref{tab:denoiser_hparams} and Table \ref{tab:regres} to facilitate reproducibility.

\begin{table}[H]
\centering
\caption{Parameters and hyperparameters for training DDPM}
\label{tab:denoiser_hparams}
\begin{tabular}{ll}
\hline
\textbf{Parameter} & \textbf{Value} \\
\hline
Input dimension ($D$)         & 44 \\
Hidden dimension ($H$)        & 256 \\
Activation                    & ReLU \\
Output dimension              & $44$ \\
Optimizer & AdamW \\
Loss function & Mean squared error (MSE) \\
\hline
\end{tabular}
\end{table}

\begin{table}[h]
\centering
\caption{ResNet regression model parameters}
\begin{tabular}{l p{0.55\columnwidth}}
\hline
\textbf{Parameter} & \textbf{Specification} \\
\hline
Hidden dimension & 256 neurons \\
Residual block & Linear(256,256) \\
Optimizer & AdamW \\
Learning rate & $5\times10^{-4}$ \\
Weight decay & $1\times10^{-4}$ \\
Batch size & 256 \\
Loss function & Mean squared error (MSE) \\
\hline
\end{tabular}
\label{tab:regres}
\end{table}

\begin{algorithm}[t]
\caption{SMC-Based Guidance }
\label{alg:smc_guidance}
\begin{algorithmic}[1]
\Require Pre-trained diffusion models $\{p_t^{\mathrm{pre}}\}_{t=T}^{0}$,
proposal policies $\{q_t\}_{t=T}^{0}$,
hyperparameter $\alpha \in \mathbb{R}$, Batch size $N$
\For{$t \in [T+1,\cdots,1]$}
  \State \textbf{Importance sampling step:} Generate $i \in [1,\cdots,N]$:
  $x_{t-1}^{[i]} \sim q_{t-1}(\cdot\mid x_t^{[i]})$
  \State Update weights:
  \Statex
  \[
    w^{[i]} \leftarrow
    \frac{
      p_{t-1}^{\mathrm{pre}}\!\left(x_{t-1}^{[i]}\mid x_t^{[i]}\right)\,
      \exp\!\left(\hat v_{t-1}\!\left(x_{t-1}^{[i]}\right)/\alpha\right)
    }{
      q_{t-1}\!\left(x_{t-1}^{[i]}\mid x_t^{[i]}\right)\,
      \exp\!\left(\hat v_t\!\left(x_t^{[i]}\right)/\alpha\right)
    }\; w^{[i]}.
  \]
  \Statex
  \State \textbf{Resampling step:} Calculate the effective sample size
  $\dfrac{\left(\sum_i w^{[i]}\right)^2}{\sum_i \left(w^{[i]}\right)^2}$.
  If it falls below 0.5N, resample by selecting new indices with replacement:
  \Statex
  \[
    \{x_{t-1}^{[i]}\}_{i=1}^N \leftarrow \{x_{t-1}^{\,s_{t-1}^{[i]}}\}_{i=1}^N,
    \qquad
    \{s_{t-1}^{[i]}\}_{i=1}^N \sim
    \mathrm{Cat}\!\left(
      \left\{
        \frac{w_{t-1}^{[i]}}{\sum_{j=1}^N w_{t-1}^{[j]}}
      \right\}_{i=1}^N
    \right).
  \]
\EndFor
\State \textbf{Output:} $\{x_0^{[i]}\}_{i=1}^N$
\end{algorithmic}
\end{algorithm}

\FloatBarrier
\input{cfd_methods}
\FloatBarrier

\section{Results}
\label{sec:results}

This section evaluates denoising accuracy and examines SMC weight concentration, providing a diagnostic of diversity in particle contributions during sampling. It then summarizes the empirical resistance screening of the generated hull forms and reports a high-fidelity hydrodynamic assessment of one screened variant and of the parametric hull representation.

\subsection{Diffusion Accuracy and Sampling Diversity}
\label{subsec:model_diagnostics}

\begin{figure}[htbp]
    \centering
    \includegraphics[width=\textwidth]{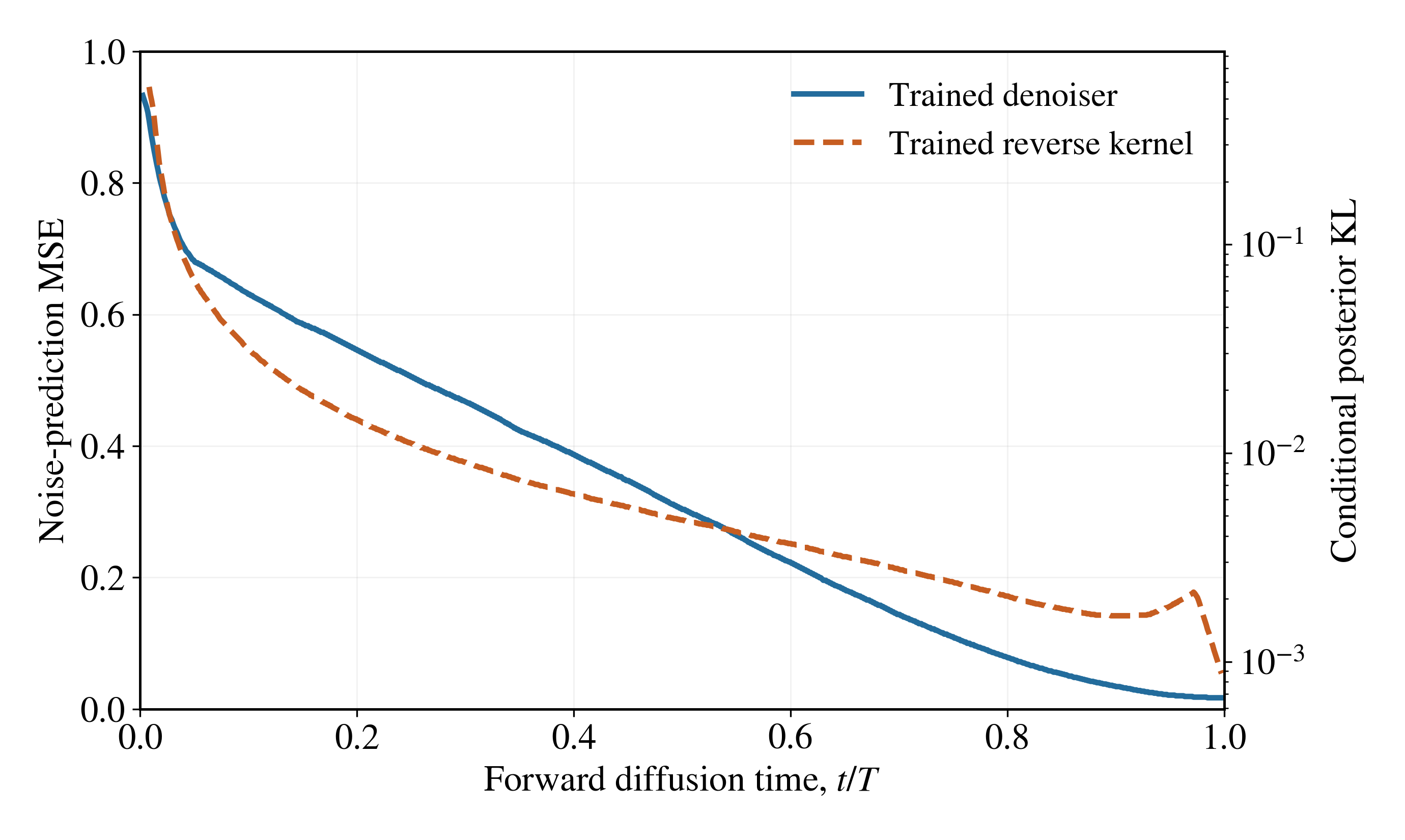}
    \caption{Denoising and reverse-step accuracy of the trained diffusion model on held-out hull designs. The left axis shows the mean squared error per dimension between the predicted and known forward-diffusion noise. The logarithmic right axis shows the conditional Gaussian posterior KL divergence, in nats per dimension, between the reference reverse-step posterior conditioned on the known clean hull and the trained reverse kernel using featurewise bounded denoised estimates. Both metrics are plotted against forward diffusion time; lower values indicate greater accuracy.}
    \label{fig:denoiser_validation}
\end{figure}

In order to assess the denoiser accuracy, we take a known, held-out hull design, add a known amount of random noise, and ask the model to predict that noise. This shows how accurately the model predicts the added noise in the forward process. Figure~\ref{fig:denoiser_validation} shows the trained model's accuracy on
held-out designs across forward diffusion time. The left axis reports
noise-prediction error, while the logarithmic right axis reports the
conditional Gaussian reverse-posterior KL divergence using the reverse
kernel with featurewise bounded denoised estimates. The KL compares the predicted reverse step distribution with the known reference hull and it shows how accurately does the model perform one reverse-diffusion step.

\FloatBarrier

\begin{figure}[htbp]
    \centering
    \includegraphics[width=0.75\textwidth]{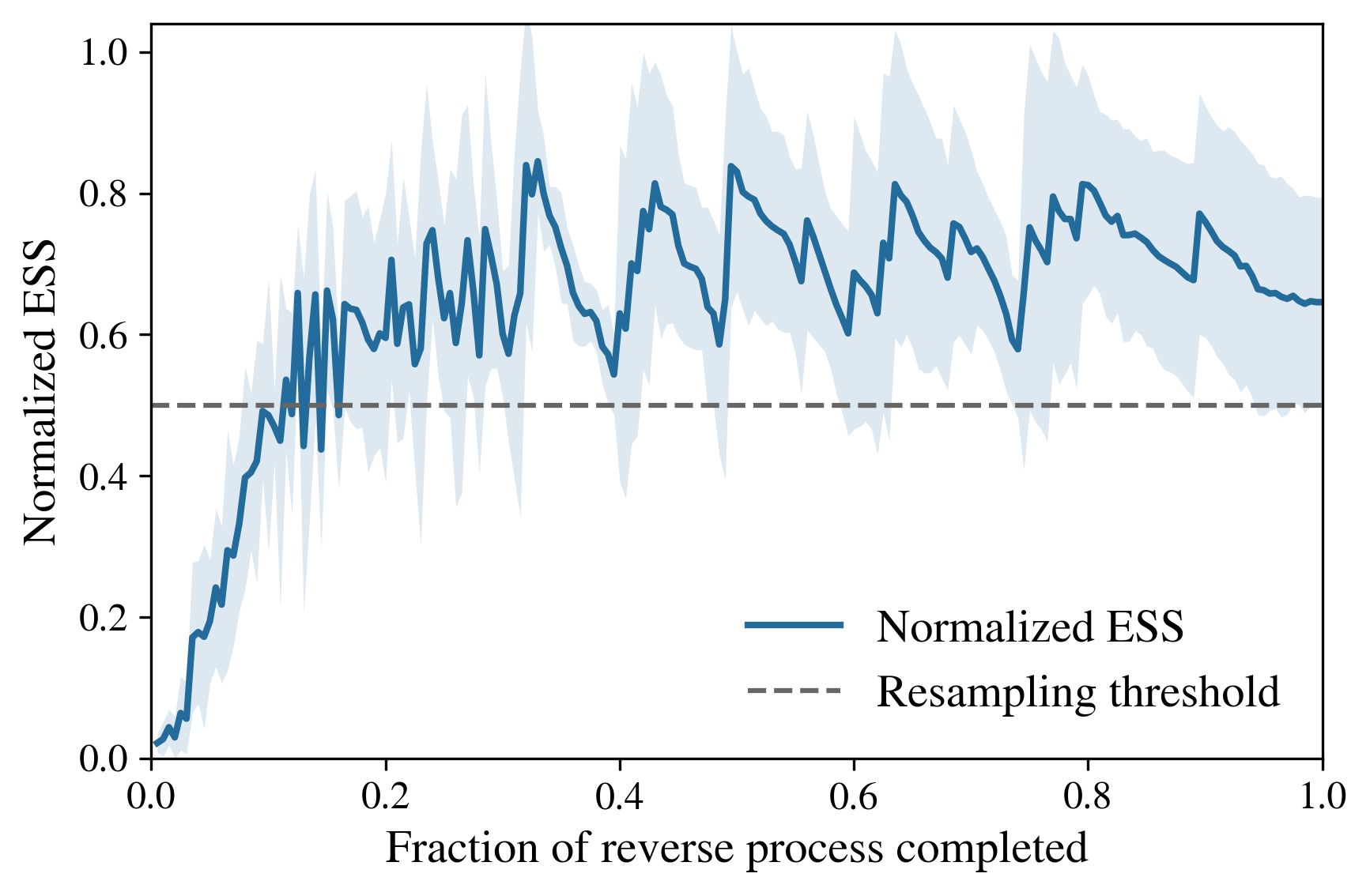}
    \caption{Normalized ESS during SMC sampling.
    The dashed line marks the resampling threshold of 0.5.
    The solid curve shows the mean and the shaded band shows an
    approximate pointwise 95\% confidence interval
    across five sampler seeds, with the trained models held fixed.
    ESS is evaluated before resampling.}
    \label{fig:smc_ess}
\end{figure}

To assess SMC sampling performance across diffusion steps, Figure~\ref{fig:smc_ess} examines the concentration of importance weights. The normalized effective
sample size, $\mathrm{ESS}/N$, is measured after each incremental
weight update and before resampling. Intermediate resampling is triggered when
$\mathrm{ESS}/N<0.5$.

\FloatBarrier

\subsection{Screening and CFD Assessment of Generated Hull Forms}
\label{subsec:empirical_screening}

This subsection summarizes the empirical screening reported in the previous
study on which the present work builds \cite{keramati2026parametric} and then assesses one of the screened variants with the CFD pipeline of Section~\ref{subsec:cfd_pipeline}.
The diffusion framework outputs candidate hull geometries as triangulated
surface meshes in STL format. Because STL is a faceted representation, the previous study applied a geometry-cleaning and surface-reconstruction step before CFD meshing. Each generated STL was imported into Rhinoceros and re-topologized
using QuadRemesh to produce a quad-dominant surface mesh with improved element regularity, converted to a SubD representation and then to NURBS surfaces,
from which a closed watertight polysurface was constructed. Geometric
integrity was verified before export, covering watertightness, consistent
surface normals and the absence of non-manifold edges and self-intersections.
The exact preprocessing operations vary with the CAD and CFD toolchain; the
requirement is a watertight, analysis-ready surface that can be robustly meshed. The generated STL of Model~1 is watertight and was meshed by the pipeline of Section~\ref{subsec:cfd_pipeline} after scaling to model size, without surface reconstruction.

Four representative variants, denoted Models~1--4, were selected from the
reward-guided sampling output and compared with the KCS baseline using the
Holtrop--Mennen method over $F_n = 0.108$ to $0.45$ at matched displacement.
All four variants returned lower predicted total resistance coefficients
than KCS over most of that range, and the largest reduction in aggregate
predicted resistance was obtained for Model~3 \cite{keramati2026parametric}.
Compared with KCS, all four variants have a reduced, effectively absent bulb
and a finer bow entrance.

These results were obtained while satisfying preliminary hydrostatic and
intact-stability screening checks, in which displacement is matched by
adjusting the draft and intact stability is evaluated through $GM$ and $GZ$
righting-arm curves against the IMO 2008 IS Code criteria. These checks
provide early-stage feasibility screening and do not constitute
class-society compliance verification. The more aggressive variants achieved
their predicted margins in part by adopting finer forms at substantially
greater draft for equal displacement. Models~3 and~4 require drafts of
approximately $17.1$ and $17.9\,\mathrm{m}$, respectively, compared with
$10.8\,\mathrm{m}$ for KCS. Such configurations would ordinarily be restricted
by route and port draft limits. The geometric characteristics of the four
variants and the influence of the parameter bounds are reported in full in
\cite{keramati2026parametric} and are not repeated here. In that study, a CFD
check of Model~1 at drafts of $0.25D$ and $0.67D$, where $D$ is the hull
depth, showed that the Holtrop--Mennen method underpredicted its resistance by
up to 40\,\% at moderate Froude number.

To compare a sampled variant with KCS under consistent conditions, Model~1
was assessed with the automated pipeline of Section~\ref{subsec:cfd_pipeline}
at the displacement of KCS, using the same mesh rules and solver settings as
the KCS computations. At this displacement, the draft of Model~1 is
$18.6\,\%$ greater than that of KCS and its waterline beam is $9.1\,\%$
greater. Its block and prismatic coefficients are $0.50$ and $0.73$, compared
with $0.65$ and $0.66$ for KCS. 
The total resistance of Model~1 differs from that of KCS by $-1.2\,\%$ at
$F_n = 0.195$, whereas it
is $44.9\,\%$ and $49.3\,\%$ higher at $F_n = 0.250$ and $0.260$,
respectively. The screening therefore selected a variant whose resistance is
comparable to that of KCS at the lower speed, but it did not capture the
increase in resistance near the design speed. Two properties of the screening are relevant to this discrepancy. The
procedure is an empirical regression evaluated on hulls whose form
coefficients may lie outside the range over which the regression was fitted,
and it aggregates over a Froude range wide enough to contain qualitatively
different regimes. The reward that guides the sampler is also a low-fidelity estimate aggregated over a wide speed range, and neither property is visible
in the values it returns.

\input{cfd_results}


\FloatBarrier
\input{cfd_discussion}
\FloatBarrier

\section{Conclusion}
\label{sec:conc}

This study combined a reward-directed diffusion framework for localized ship-hull optimization with a viscous assessment of a generated variant and of a hull refined by a naval architect toward KCS. Parametric masking and surrogate-guided sequential Monte Carlo sampling produced hull forms with lower empirically predicted calm-water resistance than KCS while restricting geometric variation to a designated region. 

An automated CFD pipeline derives mesh and solver parameters from the hull geometry and Froude number using a consistent rule set. Validated against two benchmark container ships at six speeds each, it reproduces total resistance within $+14.0\,\%$ to $-8.9\,\%$ for KCS and $-4.1\,\%$ to $-7.5\,\%$ for DTC, with errors varying monotonically with speed. The pressure converged monotonically under grid refinement, whereas the skin friction did not converge on either hull, preventing a single grid-convergence index for total resistance. These limitations support comparative assessment at matched Froude number rather than precise absolute prediction.

Applied to one of the sampled variants (Model~1), the pipeline gave a total resistance comparable to that of KCS at $F_n = 0.195$ and 45--50\,\% higher at $F_n = 0.250$ and $0.260$, an increase that the empirical screening did not capture. For the refined GenHull~E, CFD identified a speed-dependent performance reversal. Its total resistance is six to eight percent lower than that of KCS up to $F_n = 0.240$, equal within the confidence interval at $F_n = 0.250$, where an $8.07$\,N friction advantage nearly cancels an $8.49$\,N pressure penalty, and higher above that speed. The friction advantage grows smoothly with dynamic pressure, while the pressure penalty increases nearly sevenfold between $F_n = 0.240$ and $0.270$. Aggregating resistance across this crossover obscures the opposing trends.

The assessment identified limitations in geometry, reward evaluation, and stability treatment. None of the 108 bulb configurations tested produced a rounded bulbous bow; the parameterized bulb increased the resistance of GenHull~A by $27.7\,\%$ relative to KCS at $F_n = 0.152$ by viscous pressure drag. The reward omits this mechanism and uses a wetted-surface routine that undercounts sloping surfaces. GenHull~E also loses $1.04$\,m of permissible vertical center of gravity while satisfying intact-stability criteria, illustrating a margin reduction that the reward's lack of a stability term leaves unaccounted for.

These findings show that the matching of integral form coefficients and passing feasibility checks are insufficient to establish hydrodynamic performance or preserve stability margins. Continuing work will generate a CFD resistance data set of approximately 200 hulls at six speeds to train a speed-resolved reward incorporating viscous pressure and extend the assessment to further hulls sampled from the framework, before and after expert refinement. A richer geometric parameterization is also required: improved geometry and reward fidelity must advance together to generate useful features and evaluate their effects reliably.


\section*{Data and code availability}

The code for the CFD pipeline and diffusion framework, along with the dataset, will be made publicly available upon acceptance of this manuscript.


\section*{Acknowledgments}
The present study is supported by the NSERC Alliance Grant with Elomatic and InnovMarine. This research was also supported in part through computational resources and services provided by Advanced Research Computing (ARC) at the University of British Columbia and Compute Canada.

\bibliographystyle{elsarticle-num} 
\bibliography{cas-refs,cfd-validation-refs}






\end{document}

%% file: cfd_methods.tex
\subsection{CFD Evaluation of Resistance}
\label{subsec:cfd_pipeline}

The CFD results of the present study were obtained using an automated CFD pipeline developed for this purpose, in which
every mesh and solver parameter is derived by rule from the hull geometry and
the Froude number. The pipeline accepts a closed-hull STL file together with a
set of ship particulars and produces a complete OpenFOAM case, and no
parameter was adjusted manually for any of these computations.
The principal benefit of this arrangement concerns comparability more than
accuracy, since two hulls evaluated at the same Froude number are then treated
identically by construction and not according to the judgment of the analyst
preparing them. Differences between hulls may accordingly be interpreted as
differences in geometry.

\begin{figure}[htbp]
\centering
\includegraphics[width=\linewidth]{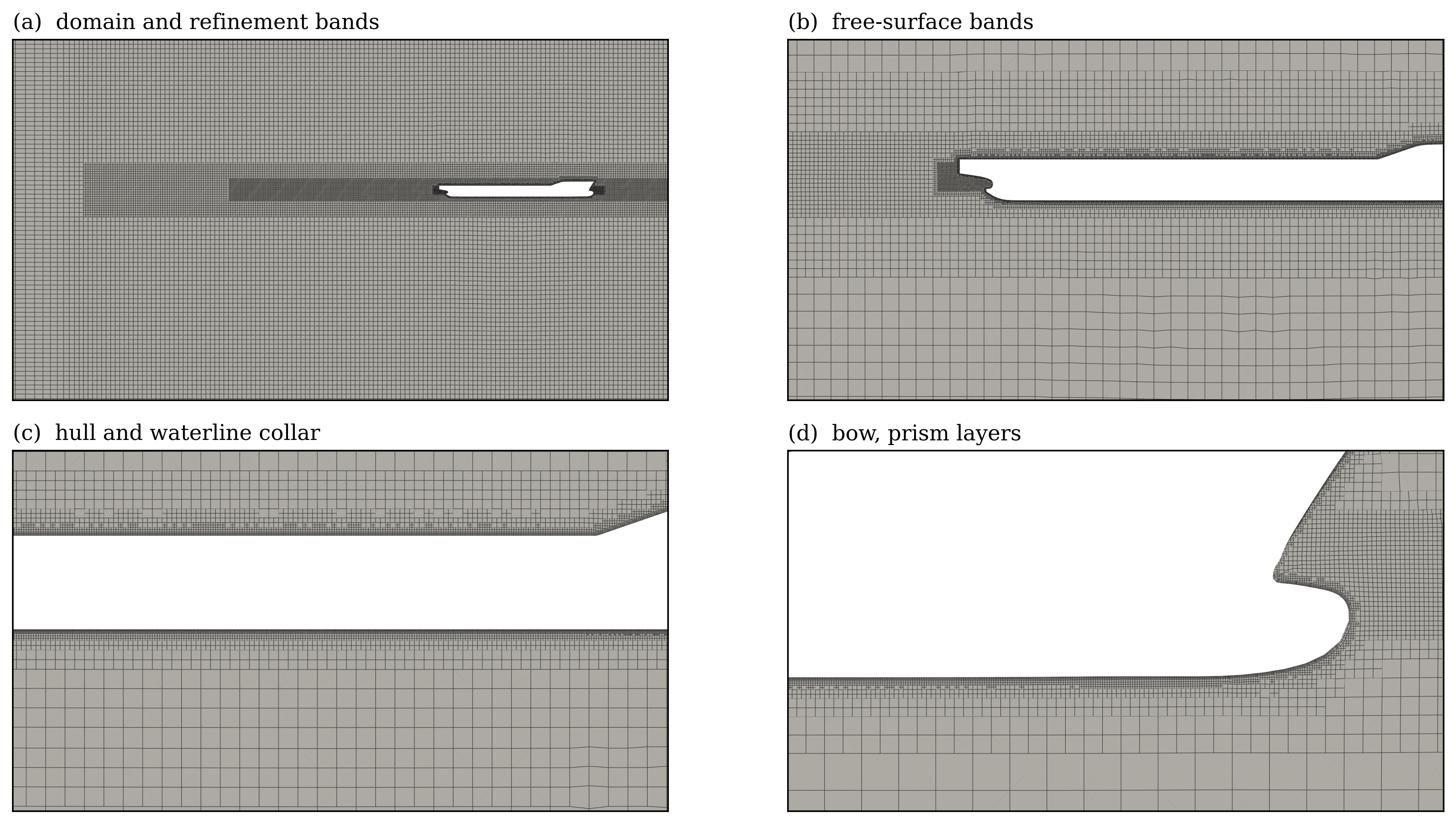}
\caption{Computational mesh generated for the KCS hull, shown on the
centerplane: (a) domain and refinement regions; (b) free-surface refinement
bands; (c) hull surface and waterline refinement; (d) prism layers at the
bow.}
\label{fig:mesh_composite}
\end{figure}

Computations use \texttt{interFoam} from OpenFOAM ESI v2412, solving the
incompressible two-phase equations with the interface captured by the
volume-of-fluid method \cite{hirt1981vof} and turbulence closed by the
$k$--$\omega$ SST model with wall functions \cite{menter2003sst}. A half
domain with a centerplane symmetry condition extends $1.5\,L_{pp}$ forward of
the forward perpendicular, $2.0\,L_{pp}$ to the side and $1.5\,L_{pp}$ below
the keel, with a Froude-scaled outlet. The background cell size is
$0.0311\,L_{pp}$, hull-surface refinement targets $y^{+}=56$, and the
free-surface band targets 64 cells across the Kelvin wavelength
$\lambda_K = 2\pi U^2/g$. The resulting meshes contain 3.52 and 9.92 million
cells for KCS and DTC respectively. Because these targets are physical rather
than integer, the same rules resolve to different refinement levels on
different hulls. A single mesh, sized at the highest Froude number of a sweep,
serves every speed in that sweep. Wave reflection from the domain boundaries
is suppressed by a sponge zone relaxing velocity and turbulence quantities
toward freestream values. Forces are averaged over at least thirty seconds of
simulated time by the method of batch means, giving a 95\,\% confidence
interval, and stationarity is verified by a half-to-half drift check.
Figure~\ref{fig:mesh_composite} shows the mesh the rules produce on the KCS
hull.

All computations are performed at fixed even keel and prescribed draft, whereas
the reference experiments are free to sink and trim. The magnitude of this
difference was measured by repeating one computation on an identical
mesh with the hull set at the experimentally measured attitude for
$F_n = 0.260$, which increased total resistance by $0.88\,\%$ and skin friction
by $1.9\,\%$. Because sinkage scales approximately with the square of speed, a
correction of this type cannot be applied uniformly across a speed range, and
the raw comparison is therefore reported throughout. The same fixed condition is applied to every hull compared in Sections~\ref{subsec:empirical_screening} and~\ref{subsec:cfd_genhull}, so that each pair is compared at the same attitude.

\subsubsection{Validation Against Benchmark Hulls}
\label{subsubsec:cfd_validation}

\begin{figure}[htbp]
\centering
\includegraphics[width=0.78\linewidth]{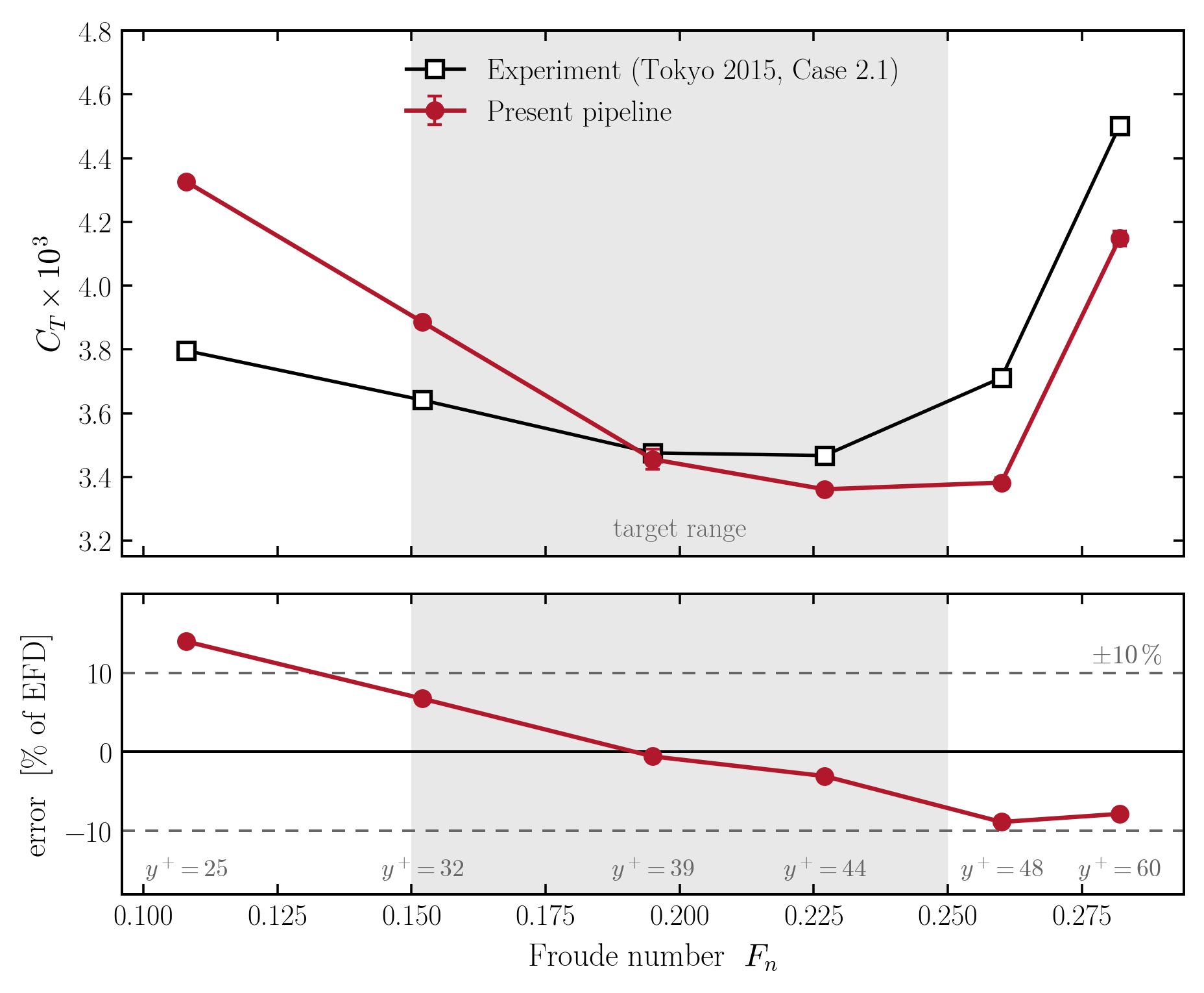}
\caption{Computed and measured total resistance coefficients for KCS (top)
and the corresponding comparison error (bottom). The shaded band indicates the Froude number range targeted for CFD-based
training data.}
\label{fig:kcs_validation}
\end{figure}

\begin{figure}[htbp]
\centering
\includegraphics[width=0.78\linewidth]{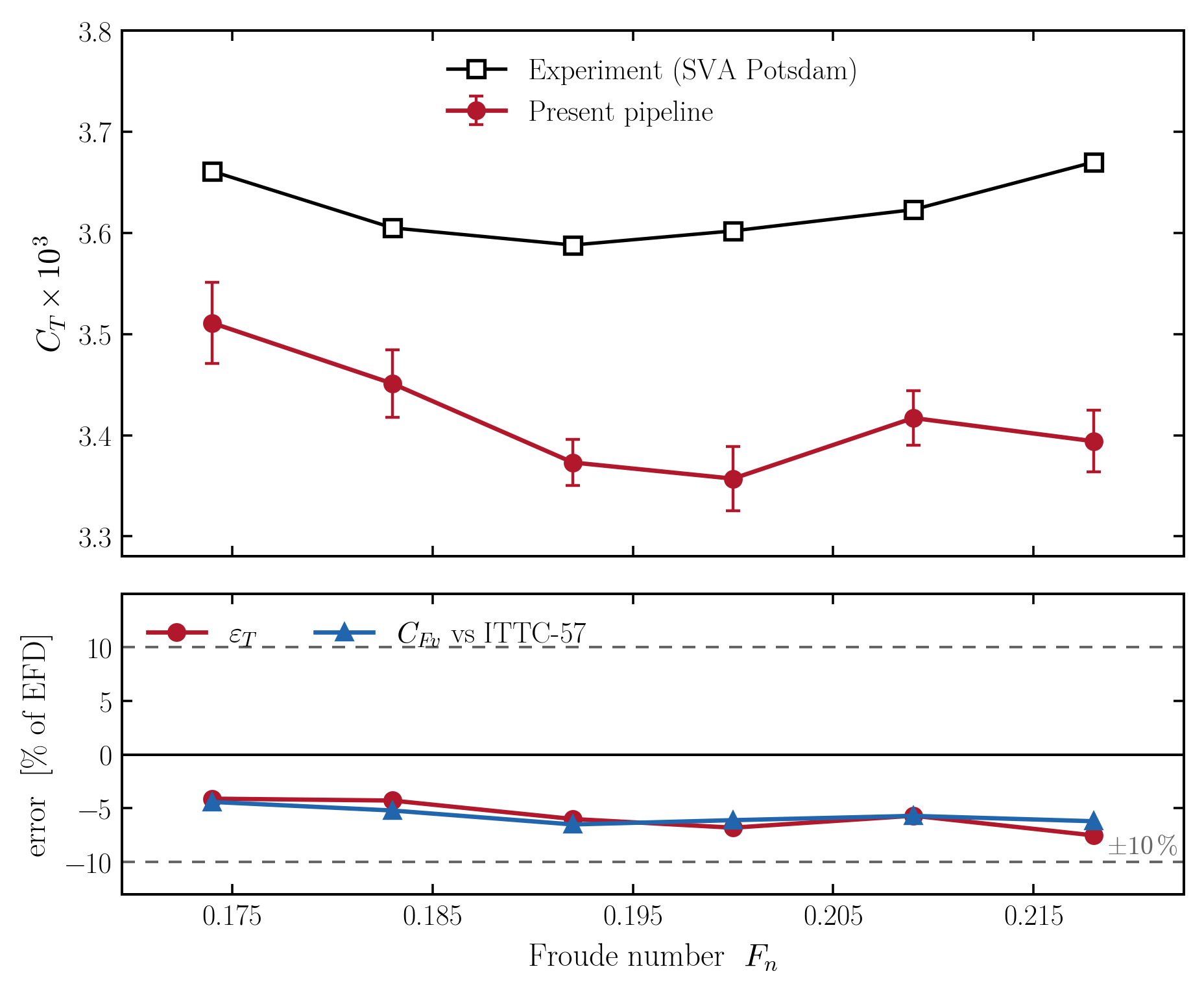}
\caption{Computed and measured total resistance coefficients for DTC (top)
and the corresponding comparison and skin-friction errors (bottom).}
\label{fig:dtc_validation}
\end{figure}

The pipeline was validated against model-scale resistance measurements for two
benchmark container ships, the KRISO Container Ship and the Duisburg Test Case
(DTC), at six Froude numbers each. The KCS results are given in
Table~\ref{tab:kcs_ladder} and Figure~\ref{fig:kcs_validation}, and the DTC
results in Table~\ref{tab:dtc_ladder} and Figure~\ref{fig:dtc_validation}.
Experimental data for KCS are taken from the Tokyo 2015 workshop, Case~2.1
\cite{tokyo2015kcs,kim2001kcs}, and for DTC from the SVA Potsdam measurements
reported by el Moctar et al.\ \cite{elmoctar2012dtc}, as tabulated by Feng
et al.\ \cite{feng2021parametric}. The KCS computations used a single mesh of
$3.52\times10^{6}$ cells, sized at $F_n = 0.282$, for all six speeds. The
tabulated $y^{+}$ is the average over the full hull patch, including dry faces
above the waterline. On the wetted surface at $F_n = 0.260$, the
distribution has a median of 98. The DTC computations used a mesh of
$9\,916\,051$ cells with $99.2\,\%$ boundary-layer coverage, generated by the
identical rule set with no hull-specific intervention, and all six cases
converged with a half-to-half drift below $0.80\,\%$. The measured DTC
resistance reaches a minimum at $F_n = 0.192$, and the computation reproduces a
minimum one speed step higher, at $F_n = 0.200$.

\begin{table}[htbp]
\centering
\caption{Validation of the CFD pipeline against KCS model-scale measurements.
$\varepsilon_T$ is the error in total resistance coefficient relative to
experiment.}
\label{tab:kcs_ladder}
\begin{tabular}{ccccccc}
\toprule
$F_n$ & $y^{+}$ & $C_T \times 10^{3}$ & 95\,\% CI & EFD $\times 10^{3}$ &
$\varepsilon_T$ [\%] & $C_{Fv}$ vs ITTC-57 [\%] \\
\midrule
0.108 & 25 & 4.326 & $\pm0.08\,\%$ & 3.796 & $+14.0$ & $+6.3$ \\
0.152 & 32 & 3.886 & $\pm0.40\,\%$ & 3.641 & $+6.7$  & $+2.9$ \\
0.195 & 39 & 3.455 & $\pm0.90\,\%$ & 3.475 & $-0.6$  & $-2.1$ \\
0.227 & 44 & 3.361 & $\pm0.24\,\%$ & 3.467 & $-3.1$  & $-4.5$ \\
0.260 & 48 & 3.382 & $\pm0.46\,\%$ & 3.711 & $-8.9$  & $-6.3$ \\
0.282 & 60 & 4.148 & $\pm0.55\,\%$ & 4.501 & $-7.8$  & $-7.8$ \\
\bottomrule
\end{tabular}
\end{table}

\begin{table}[htbp]
\centering
\caption{Validation of the CFD pipeline against DTC model-scale measurements,
using the same rule set as for KCS.}
\label{tab:dtc_ladder}
\begin{tabular}{ccccccc}
\toprule
$F_n$ & $Re$ & $C_T \times 10^{3}$ & 95\,\% CI & EFD $\times 10^{3}$ &
$\varepsilon_T$ [\%] & $C_{Fv}$ vs ITTC-57 [\%] \\
\midrule
0.174 & $7.30\times10^{6}$ & 3.511 & $\pm1.15\,\%$ & 3.661 & $-4.1$ & $-4.4$ \\
0.183 & $7.68\times10^{6}$ & 3.451 & $\pm0.96\,\%$ & 3.605 & $-4.3$ & $-5.2$ \\
0.192 & $8.06\times10^{6}$ & 3.373 & $\pm0.68\,\%$ & 3.588 & $-6.0$ & $-6.5$ \\
0.200 & $8.40\times10^{6}$ & 3.357 & $\pm0.95\,\%$ & 3.602 & $-6.8$ & $-6.1$ \\
0.209 & $8.77\times10^{6}$ & 3.417 & $\pm0.78\,\%$ & 3.623 & $-5.7$ & $-5.7$ \\
0.218 & $9.15\times10^{6}$ & 3.394 & $\pm0.90\,\%$ & 3.670 & $-7.5$ & $-6.2$ \\
\bottomrule
\end{tabular}
\end{table}

The KCS points stay between $+14.0\,\%$ and $-8.9\,\%$ of the measurement,
and between $+6.7\,\%$ and $-3.1\,\%$ over the range $F_n = 0.152$ to $0.227$
targeted for the CFD-based training data, while DTC points vary between
$-4.1\,\%$ and $-7.5\,\%$. Both hulls show skin friction trending with speed, at $+6.3\,\%$
to $-7.8\,\%$ on KCS and $-4.4\,\%$ to $-6.2\,\%$ on DTC.

The comparison error varies monotonically with Froude number, and this
behavior follows from the meshing rules. Because a single mesh serves
an entire sweep, the near-wall spacing remains fixed while the speed varies,
so that the skin-friction error follows the wall resolution. The free-surface
band is similarly sized for the highest speed of the sweep, so that the
wave-resistance error follows the number of cells resolving the Kelvin
wavelength, which falls from 64 at the top of the KCS range to approximately
nine at the bottom. The comparison error therefore constitutes a bounded
envelope with an identified cause, and may be taken into account when the
comparisons presented below are interpreted. The lowest KCS speed,
$F_n = 0.108$, lies outside the range used for any comparison in this paper
and has been retained because it provides the coarse-resolution end of this
explanation.

Two limitations of the validation are noted. A grid-convergence study was
carried out for each hull on three geometrically similar grids, with the
background cell size scaled by a factor of approximately $\sqrt{2}$, and the
convergence ratio was evaluated separately for the pressure and viscous
components of the resistance. The viscous component did not converge under
refinement on either hull, with convergence ratios of 2.74 for KCS and 1.38
for DTC, whereas the pressure component converged monotonically with ratios of
0.68 and 0.52 respectively. The behavior of the total resistance under
refinement follows that of the dominant component, and is divergent for KCS
and convergent for DTC. A grid-convergence index for total resistance could
therefore not be established and is not reported. The second limitation is
that the comparison error at the KCS design speed exceeds the validation
uncertainty, which indicates the presence of modeling error according to the
ITTC verification and validation procedure \cite{ittc2021verification}. The
residual is consistent with the behavior of the wall functions at the
achieved $y^{+}$, although this was not confirmed by a wall-resolved
computation.

The pipeline is therefore suitable for comparative ranking of hulls at matched
Froude number under a single rule set, in which the systematic error is common
to both hulls and largely cancels in the difference, and not for absolute
prediction of resistance. The degree of cancellation is bounded by the present
validation data. At the closest matched speeds the comparison errors of the
two benchmark hulls differ by four to five percentage points, so that
cancellation cannot be assumed between dissimilar hull forms. For two hulls of
matched length, draft and displacement, such as KCS and GenHull~E (Section~\ref{subsec:cfd_genhull}), the errors are expected to be more nearly equal and a greater degree of cancellation may be anticipated.

%% file: cfd_results.tex
\subsection{CFD Assessment of a Refined Hull}
\label{subsec:cfd_genhull}

The discrepancy between the empirical and viscous assessments of Model~1, first observed in the earlier study \cite{keramati2026parametric}, led to an expert examination of the parametric representation of Section~\ref{sec:methodology}. Because the sampler can propose only hulls that this representation is able to express, its capability bounds what the framework can achieve. In this section, a hull refined by a naval architect toward KCS is evaluated with the same pipeline (Section~\ref{subsec:cfd_pipeline}) to examine what the representation produces near the parent hull and how the resistance of that hull varies across the speed range. The hull considered, designated GenHull~E, is distinct from Models~1--4. It was constructed by the authors as an expert check of the representation and was fitted to KCS by constrained least squares, starting from a KCS-like parameter set based on the
principal dimensions of KCS. In the first of two fitting stages, nine integral form coefficients were matched, with length, draft and displacement imposed as hard constraints and the remaining quantities, including beam, midship area, longitudinal center of buoyancy and waterplane coefficient, entering as soft residuals. In the second, the station half-breadths of KCS were matched
together with waterline length, displacement and beam, with the aft stations weighted most heavily and the bulb parameters disabled. No resistance or stability term entered the fit at either stage. GenHull~E is therefore the closest form to KCS that the representation was found to express, and its behavior characterizes the representation within which the generative framework operates, not the output of the sampler. GenHull~A, an earlier stage of the same sequence that retained the
parameterized bulb, is discussed in Section~\ref{subsec:representation_limits}.

\subsubsection{Geometry and Hydrostatics}
\label{subsubsec:genhull_geometry}

\begin{figure}[htbp]
\centering
\includegraphics[width=\linewidth]{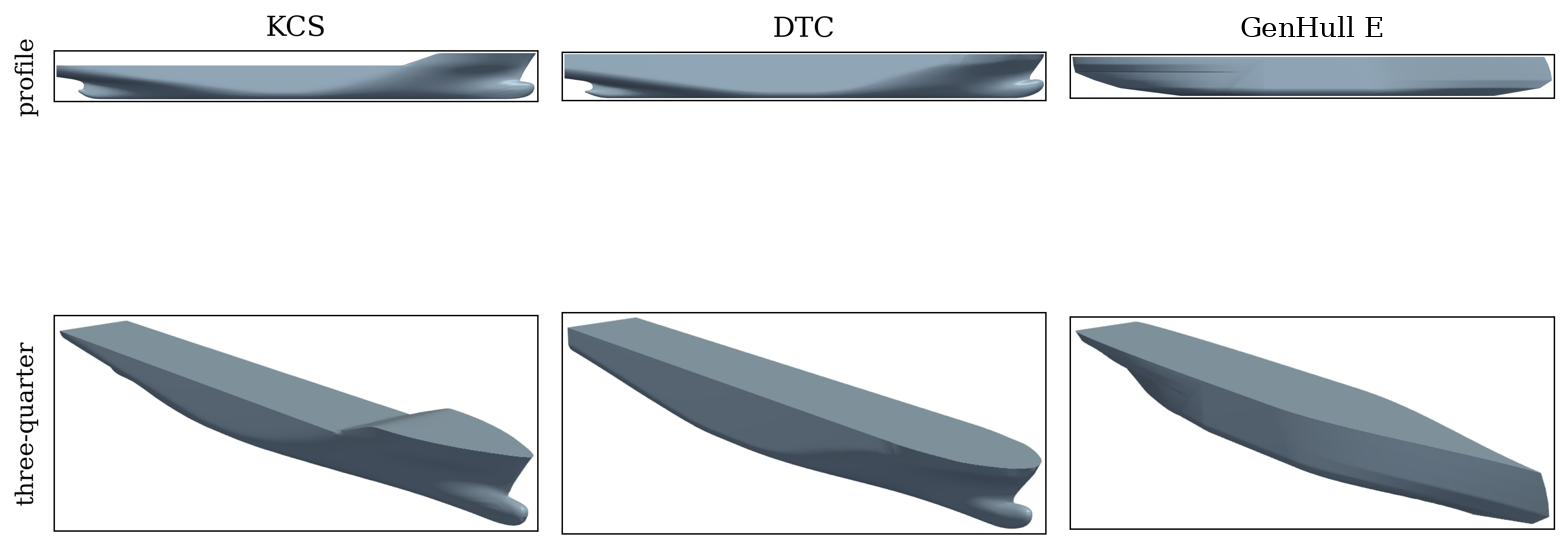}
\caption{Profile (top) and three-quarter (bottom) views of the KCS, DTC and
GenHull~E hull forms.}
\label{fig:hull_composite}
\end{figure}

Figure~\ref{fig:hull_composite} shows GenHull~E together with the two
benchmark container ships used to validate the evaluation pipeline, and
Table~\ref{tab:hull_particulars} lists the principal particulars and form
coefficients of GenHull~E and KCS. All values were computed from the as-meshed
geometry at the design draft. The coefficients are referred to $L_{pp}$, which
is the basis for the Froude and Reynolds numbers and is common to both hulls
by construction. The waterline length is reported separately, because the KCS
bulb protrudes forward of the forward perpendicular and the measured waterline
consequently exceeds $L_{pp}$ by $0.95\,\%$. The longitudinal center of
buoyancy is measured identically on both hulls as a percentage of waterline
length aft of the waterline mid-point; the published KCS value of
$+1.48\,\%L_{pp}$ is referred to the forward perpendicular, which the bulb
displaces.

\begin{table}[htbp]
\centering
\caption{Principal particulars and form coefficients of KCS and GenHull~E at
full scale, computed from the as-meshed geometry at design draft.}
\label{tab:hull_particulars}
\begin{tabular}{lrrr}
\toprule
Quantity & KCS & GenHull~E & Difference \\
\midrule
$L_{pp}$ [m]              & 230.00  & 230.00  & constrained \\
$L_{wl}$ [m]              & 232.19  & 230.30  & $-0.81\,\%$ \\
$B_{wl}$ [m]              & 32.22   & 34.46   & $+6.96\,\%$ \\
$T$ [m]                   & 10.80   & 10.80   & $-0.01\,\%$ \\
$\nabla$ [m$^{3}$]        & 51\,765 & 52\,042 & $+0.53\,\%$ \\
$S$ [m$^{2}$]             & 9\,462  & 8\,802  & $-6.97\,\%$ \\
$A_{wp}$ [m$^{2}$]        & 6\,135  & 5\,466  & $-10.91\,\%$ \\
$A_{m}$ [m$^{2}$]         & 341.1   & 360.9   & $+5.80\,\%$ \\
$C_B$                     & 0.6468  & 0.6080  & $-6.00\,\%$ \\
$C_M$                     & 0.9802  & 0.9697  & $-1.08\,\%$ \\
$C_P$                     & 0.6599  & 0.6270  & $-4.98\,\%$ \\
$C_{WP}$                  & 0.8279  & 0.6896  & $-16.71\,\%$ \\
$S/\nabla^{2/3}$          & 6.812   & 6.315   & $-7.30\,\%$ \\
LCB [\%$L$ aft]           & $+1.00$ & $+0.85$ & 0.15 forward \\
$BM_T$ [m]                & 9.00    & 7.84    & $-12.93\,\%$ \\
$KM_T$ [m]                & 14.93   & 13.64   & $-8.62\,\%$ \\
Bulbous bow               & present & absent  & --- \\
\bottomrule
\end{tabular}
\end{table}

At matched length, draft and displacement, GenHull~E has $6.97\,\%$ less wetted surface than KCS, $6.96\,\%$ greater beam, and a finer entry. These differences are consequences of the representation and not of any objective, since the fit contained no resistance term. With the bulb disabled and the hull fitted to the station half-breadths of KCS, the closest
form that the representation produces at matched length, draft and displacement is beamier and finer at the ends than KCS and carries less wetted surface. GenHull~E has no bulbous bow, whereas both benchmark hulls possess one, because the representation cannot express a faired bulb, as examined in
Section~\ref{subsec:representation_limits}.

\begin{figure}[htbp]
\centering
\includegraphics[width=\linewidth]{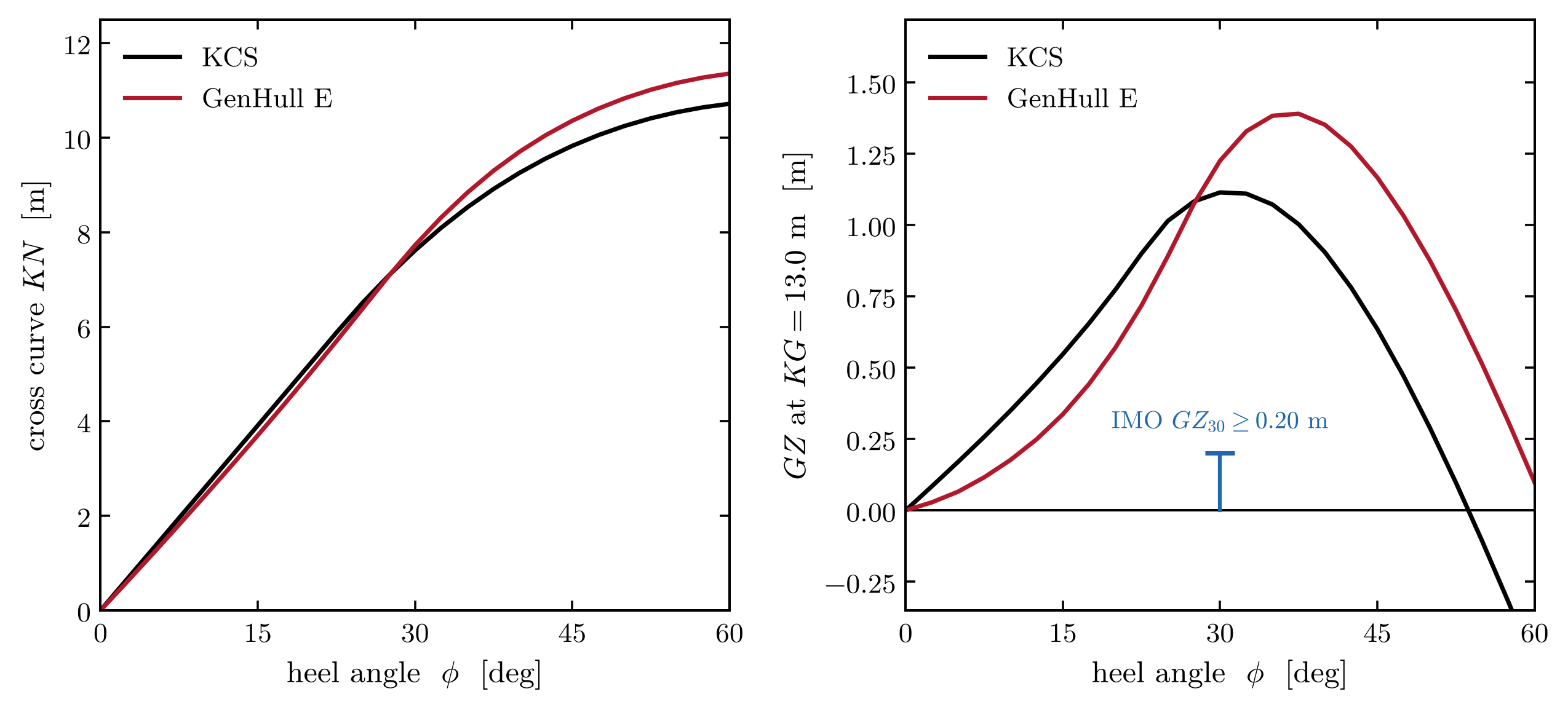}
\caption{Cross curves of stability (left) and righting-arm curves at
$KG = 13.0$\,m (right) for KCS and GenHull~E at design draft.}
\label{fig:stability}
\end{figure}

The finer waterplane has a consequence for transverse stability that is not
apparent from the form coefficients. Although the beam increases by
$6.96\,\%$, the waterplane inertia decreases by $12.46\,\%$, because a fine
entry places less waterplane area at large distance from the centerline. The
metacentric radius $BM_T$ consequently decreases by $12.93\,\%$ and the height
of the transverse metacentre $KM_T$ by $1.29$\,m. Figure~\ref{fig:stability}
shows the cross curves of stability computed for both hulls at constant
displacement, with both hulls truncated to a common depth of 19.0\,m so that
the freeboard treatment is identical. The two curves cross near $28^\circ$.
GenHull~E is less stiff than KCS near the upright condition but develops a
larger maximum righting lever at a greater angle of heel, reaching
$GZ_{max} = 1.39$\,m at $37.5^\circ$ against $1.11$\,m at $30^\circ$ for KCS
when both hulls are evaluated at $KG = 13.0$\,m.

Expressed as the highest center of gravity satisfying the general
intact-stability criteria of the 2008 IS Code, the limiting $KG$ is
$14.53$\,m for KCS and $13.49$\,m for GenHull~E, a reduction of $1.04$\,m.
The governing criterion differs between the two hulls. KCS is limited by the
area under the righting-arm curve between $30^\circ$ and $40^\circ$, whereas
GenHull~E is limited by the requirement for a minimum initial metacentric
height. Where an operational metacentric height requirement applies equally to
both hulls, as is usual for a container ship, the reduction corresponds to the
full $1.29$\,m difference in $KM_T$. All stability quantities are reported at
full scale, since the criteria of the IS Code are defined at that scale. The
values were computed at upright trim on a bare hull with no downflooding angle
defined, and the severe wind and rolling criterion was not evaluated. They
therefore constitute a comparison of hull form and not an assessment of
compliance.

\subsubsection{Resistance Comparison with KCS}
\label{subsubsec:genhull_resistance}

\begin{figure}[htbp]
\centering
\includegraphics[width=\linewidth]{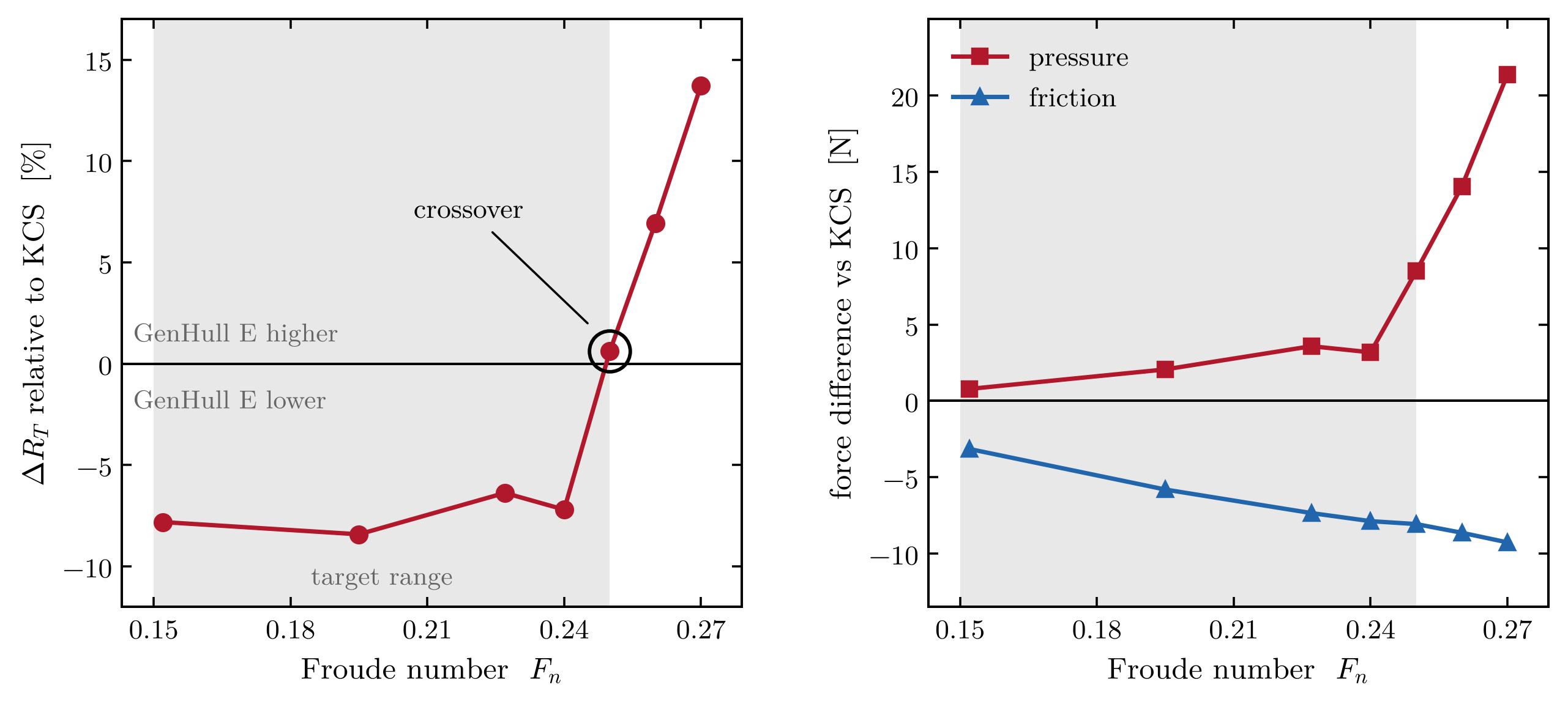}
\caption{Total resistance of GenHull~E relative to KCS (left) and the pressure
and friction components of the difference (right). The shaded band indicates the Froude number range targeted for CFD-based
training data.}
\label{fig:genhull_crossover}
\end{figure}

GenHull~E was compared with KCS at seven Froude numbers, with length, draft
and displacement as the matched basis and with beam and wetted surface free.
Both hulls were computed on the same machine under the identical rule set, so
that no difference in numerical setup enters the comparison, and the forces on
both were averaged over matched windows of the force histories. The results
are given in Table~\ref{tab:genhull_resistance} and
Figure~\ref{fig:genhull_crossover}. The final column of the table gives the
resistance difference divided by the combined 95\,\% confidence half-width of
the pair; it is a measure of separation and not a standard-deviation
significance level, and it is not reported at the two lowest speeds, for which
the confidence intervals were not recorded.

\begin{table}[htbp]
\centering
\caption{Total resistance of GenHull~E and KCS at model scale $1{:}31.6$, with
the pressure and friction components of the difference.}
\label{tab:genhull_resistance}
\small
\setlength{\tabcolsep}{5pt}
\begin{tabular}{ccccccc}
\toprule
$F_n$ & $R_T$ KCS & $R_T$ GenHull~E & $\Delta R_T$ & $\Delta$ pressure &
$\Delta$ friction & separation \\
 & [N] & [N] & [\%] & [N] & [N] & \\
\midrule
0.152 & 30.418 & 28.042  & $-7.81$  & $+0.78$  & $-3.15$ & --- \\
0.195 & 44.596 & 40.843  & $-8.42$  & $+2.06$  & $-5.81$ & --- \\
0.227 & 59.147 & 55.377  & $-6.37$  & $+3.58$  & $-7.35$ & 2.1 \\
0.240 & 65.366 & 60.661  & $-7.20$  & $+3.18$  & $-7.88$ & 1.9 \\
0.250 & 70.053 & 70.480  & $+0.61$  & $+8.49$  & $-8.07$ & 0.2 \\
0.260 & 77.716 & 83.114  & $+6.95$  & $+14.04$ & $-8.64$ & 3.4 \\
0.270 & 88.128 & 100.231 & $+13.73$ & $+21.37$ & $-9.27$ & 8.4 \\
\bottomrule
\end{tabular}
\end{table}

GenHull~E has six to eight percent lower total resistance than KCS from
$F_n = 0.152$ to $0.240$ and higher total resistance above that range. The
crossover occurs at $F_n = 0.250$, where the difference is $+0.61\,\%$, or
$0.43$\,N on $70$\,N. This difference is two tenths of the combined confidence
half-width and is therefore indistinguishable from zero. When evaluated on
four successive averaging windows as the run lengthened, the difference at
this speed read $+0.27$, $+0.02$, $+0.22$ and $+0.61\,\%$, all within half a
percent of zero. The crossover is thus located at a computed speed and not
interpolated between distant points, and the four window lengths place it at
$F_n = 0.250$ to within the repeatability of the method.

The component columns of Table~\ref{tab:genhull_resistance} identify the
mechanism responsible. GenHull~E has lower friction than KCS at every speed,
and the difference increases smoothly with dynamic pressure from $-3.15$\,N at
$F_n = 0.152$ to $-9.27$\,N at $0.270$. The magnitude of this difference
exceeds the value that the reduced wetted surface alone would produce, since
at the crossover the friction difference amounts to $14.3\,\%$ of the KCS
value against a $6.97\,\%$ difference in wetted area. This excess is examined
in Section~\ref{subsubsec:genhull_limitations}. The pressure difference
remains below $4$\,N up to $F_n = 0.240$ and then increases sharply, reaching $8.5$, $14.0$ and $21.4$\,N at $0.250$, $0.260$ and $0.270$ respectively. Between $F_n = 0.240$ and $0.270$ it increases nearly sevenfold, whereas the dynamic pressure increases by $27\,\%$. The rate of this increase identifies the pressure difference as
wave-making resistance and not as viscous pressure resistance. The two
components are measured separately and have opposite speed dependence, and the
crossover occurs at the speed where the pressure penalty exceeds the friction
advantage.

The matched-displacement basis has a further consequence for the manner in
which the comparison is reported. Because the wetted surfaces differ by
approximately $7\,\%$, the total resistance coefficient and the total
resistance force do not rank the hulls by the same margin, and at
$F_n = 0.270$ a coefficient difference of $22\,\%$ corresponds to a force
difference of $13.7\,\%$. Forces are therefore reported in this section, and
coefficients are reserved for the comparison with measurement in
Section~\ref{subsubsec:cfd_validation}.

\subsubsection{Free-Surface and Hull Pressure Fields}
\label{subsubsec:genhull_flow}

\begin{figure}[htbp]
\centering
\includegraphics[width=\linewidth]{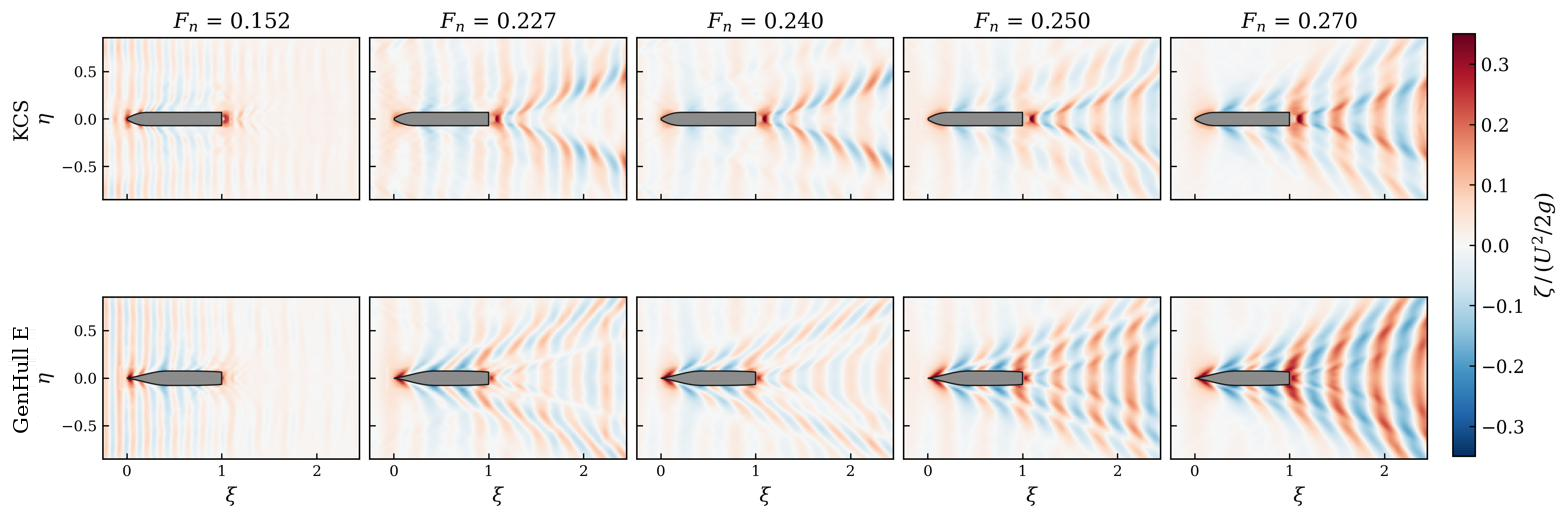}
\caption{Free-surface elevation, normalized by $U^{2}/2g$, for KCS (top) and
GenHull~E (bottom) at five Froude numbers.}
\label{fig:wave_crossover}
\end{figure}

Figure~\ref{fig:wave_crossover} shows the computed free-surface elevation for
both hulls at five Froude numbers and indicates the physical origin of the
crossover. At $F_n = 0.152$ the Kelvin wavelength is $0.145\,L_{pp}$ and the
two wave systems are nearly indistinguishable, so that wave-making forms a
small fraction of the total resistance and the difference in wetted surface
determines the comparison. At $F_n = 0.270$ the wavelength has increased to
$0.458\,L_{pp}$, both wave systems are strong, and the system generated by
GenHull~E is visibly the stronger. The elevation is normalized by $U^{2}/2g$
because wave amplitude scales with the square of speed; on a common
dimensional scale the low-speed panels would be nearly blank and differences
between speeds would obscure differences between hulls.

\begin{figure}[htbp]
\centering
\includegraphics[width=\linewidth]{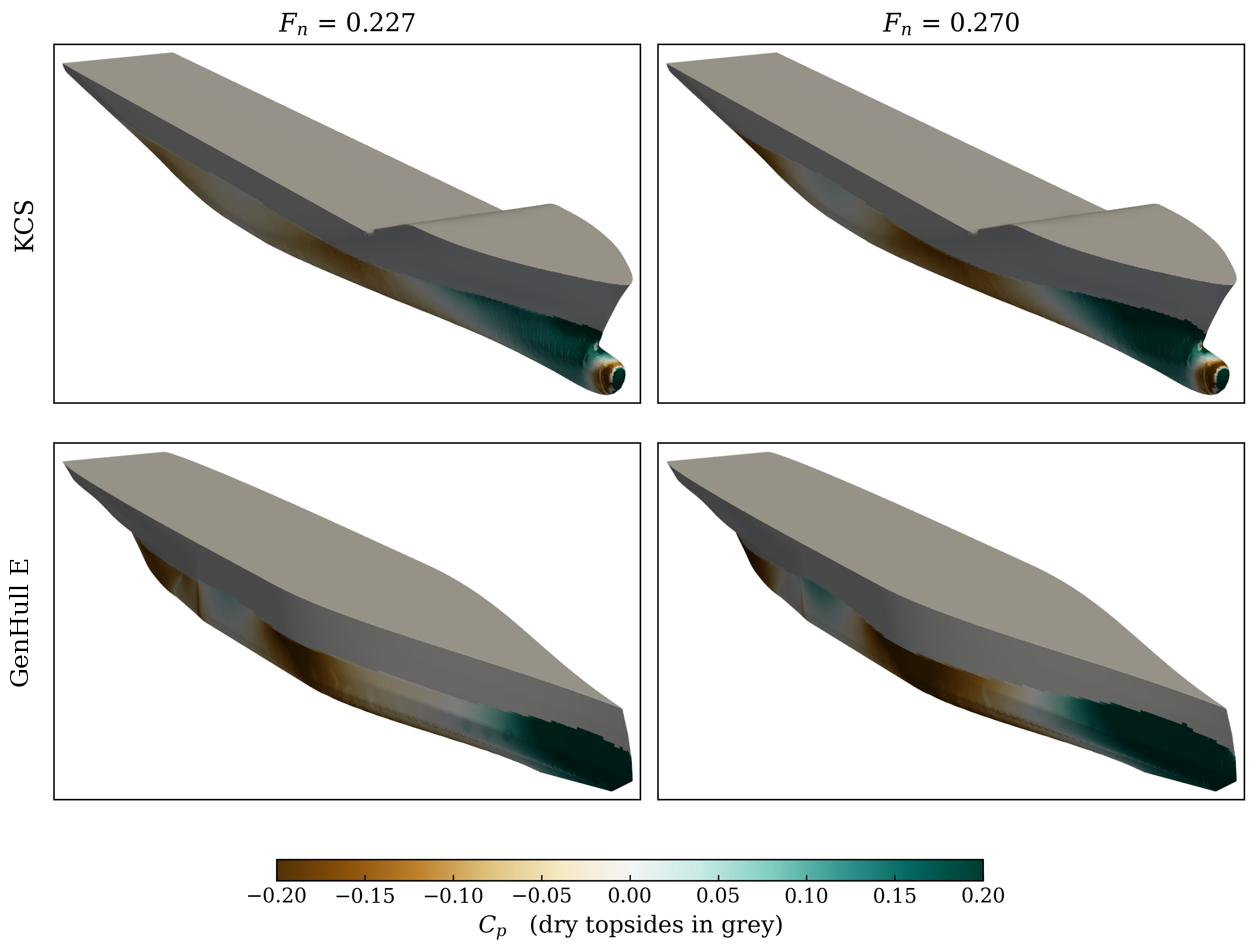}
\caption{Pressure coefficient on the wetted forebody of KCS (top) and
GenHull~E (bottom) at $F_n = 0.227$ and $0.270$.}
\label{fig:cp_comparison}
\end{figure}

Figure~\ref{fig:cp_comparison} shows the pressure coefficient on the wetted
forebody of both hulls at Froude numbers on either side of the crossover. KCS
exhibits a concentrated stagnation region ahead of the stem at the bulb, which
is the wave-cancellation mechanism that a bulbous bow provides near the design
speed. GenHull~E exhibits its stagnation region on the stem and develops
stronger forebody suction at the higher speed, with no canceling wave system
ahead of the hull.

\subsubsection{Limitations of the Comparison}
\label{subsubsec:genhull_limitations}

Two aspects of the computed viscous force qualify the comparison. The first concerns boundary-layer coverage, which lies between 87 and $92\,\%$ on GenHull~E against $98\,\%$ on KCS. A direct measurement on the DTC hull, in which meshes achieving $75.8\,\%$ and $99.2\,\%$ coverage were compared at
the same condition, established that 23.4 percentage points of coverage correspond to $3.7\,\%$ of skin friction. Scaling this result linearly, the present shortfall of approximately 11 points corresponds to $1.7\,\%$ of skin
friction, or approximately $1.4\,\%$ of total resistance, in favor of GenHull~E. A contribution of this size is attributable to discretization and should be discounted from the measured advantage.

The second concerns the remaining discrepancy, which is larger and has a different origin. At the crossover the viscous force on GenHull~E lies $13.4\,\%$ below the ITTC-57 flat-plate value, against $6.1\,\%$ for KCS, a gap
of 7.3 percentage points of which coverage accounts for fewer than two. A station-wise diagnostic located the balance of this gap at the bow and the forward shoulder, which is consistent with local separation on a forebody that
has no bulb and carries excess beam at the shoulder. Separation of this kind reduces wall shear and increases pressure drag simultaneously, so that its effect is not confined to a reduction in friction, and the associated increase
in pressure drag is included in the pressure column of
Table~\ref{tab:genhull_resistance}. The present meshes do not permit the physical and numerical contributions to be distinguished, and no correction has been applied. When only the coverage contribution is discounted, the
resistance advantage attributable to the hull form up to $F_n = 0.240$ is five to seven percent, against the six to eight percent measured.

No grid-convergence study was performed on GenHull~E. The meshing rules applied to both hulls are identical to those used for the validation cases, so that the discretization error is largely common to the two members of each pair, although this argument is weakest where the flow physics differs
most, which is above the crossover. Both hulls were computed at fixed even keel, consistent with the validation cases.

%% file: cfd_discussion.tex
\section{Discussion}
\label{sec:discussion}

\subsection{Influence of Froude Number on the Resistance Comparison}
\label{subsec:speed_dependence}

The comparison of Section~\ref{subsubsec:genhull_resistance} shows that the sign of the resistance difference between GenHull~E and KCS depends on the Froude number at which the comparison is made. At $F_n = 0.195$ GenHull~E has $8.4\,\%$ lower total resistance than the benchmark, whereas at $F_n = 0.270$ its resistance is $13.7\,\%$ higher. Model~1, the sampled variant of Section~\ref{subsec:empirical_screening}, shows a similar dependence on speed with a larger penalty at the higher speeds. Its resistance is comparable to that of KCS at $F_n = 0.195$ but 45--50\,\% higher at $F_n = 0.250$ and $0.260$, whereas the empirical screening predicted a lower resistance coefficient than that of KCS near $F_n = 0.26$ \cite{keramati2026parametric}. 
For both hulls, an assessment made at a single speed would therefore not hold across the speed range.

For GenHull~E, the reversal arises from the two mechanisms identified in Section~\ref{subsubsec:genhull_resistance}, in which the friction advantage of the smaller wetted surface grows smoothly with dynamic pressure while the wave-making penalty of the bulbless forebody grows much faster above $F_n = 0.240$. An objective aggregated over a Froude range
that spans such a reversal assigns a single value to two opposing regimes, and the performance at any individual speed may differ substantially from that value. Locating the speed at which an advantage ends requires evaluation at a sufficient number of conditions to bracket the crossover, and in the present case seven paired computations were required to establish it. An automated evaluation pipeline makes computations
of this number affordable for a single hull comparison.

The result also bears on the selection of the evaluation range within the generative loop. The reward of Section~\ref{sec:methodology}, following 
\cite{keramati2026reward}, aggregates resistance from $F_n = 0.10$ to $0.45$, a range that spans slow-steaming conditions through the resistance hump characteristic of container-ship forms. Because resistance grows rapidly with speed, an aggregate of this kind is dominated by the highest speeds, which here
lie above the crossover. The operational range of a given vessel is usually narrower. For a vessel operating predominantly below $F_n = 0.2$, weighting the objective toward the low-speed end would favor forms resembling GenHull~E. For a container ship operating between $F_n = 0.25$ and $0.30$, the CFD results indicate that the same form would be disadvantageous, and that the transition between the two regimes is sufficiently sharp that an objective averaged across it would not represent either regime adequately. Selection of the Froude number range and its associated weighting is therefore an essential part of the problem definition.

\subsection{Bulbous Bow Representation in the Parametric Hull Model}
\label{subsec:representation_limits}

GenHull~E has no bulbous bow, and this characteristic was found to originate in the parametric hull representation and not in the fit. A sweep of 108 parameter configurations was carried out over the bulb-related parameters of the representation, and every configuration produced a constant-width
cylindrical forebody, since the bulb profile blends against a section family that cannot close into a rounded forward body. One of the nominal bulb parameters was found to have no effect on the generated geometry.

The consequence of this limitation was quantified at $F_n = 0.152$ under the same meshing rules. GenHull~A, an earlier stage of the fit that retained the parameterized bulb while matching KCS on the nine integral form coefficients, returned $27.7\,\%$ higher total resistance than
KCS, whereas GenHull~E, with the bulb disabled, returned $7.8\,\%$ lower resistance at the same speed. The excess resistance of GenHull~A was traced to viscous pressure drag at its cylindrical forebody.

This observation has a general implication for generative hull models. Matching integral form coefficients does not constrain the local geometry that governs resistance, and a parameterization may be incapable of representing a
feature without any indication of this from the coefficients on which it was fitted. Verification of a generative geometry model against integral coefficients alone is therefore insufficient, and verification against local geometry or against a resolved resistance computation is required in addition.

\subsection{Transverse Stability and the Objective Function}
\label{subsec:constraint_objective}

The stability result of Section~\ref{subsubsec:genhull_geometry} illustrates a related limitation. The generative framework screens candidate geometries against the general intact-stability criteria of the 2008 IS Code, and GenHull~E satisfies those criteria for a vertical center of gravity up to
$13.49$\,m. The fit that produced GenHull~E constrained integral form coefficients and station shapes, none of which is a measure of stability, and the reward of the generative framework likewise contains no stability term. The resulting hull has a limiting vertical center of gravity $1.04$\,m lower
than that of the benchmark, or $1.29$\,m lower where an operational metacentric height requirement governs in place of the regulatory minimum, corresponding to seven to nine percent of a limiting $KG$ of approximately 14.5\,m.

The cause of this reduction is not apparent from the form coefficients. GenHull~E is $6.96\,\%$ beamier than KCS, which would ordinarily increase the metacentric radius, but its waterplane coefficient is $16.7\,\%$ lower and the waterplane inertia is $12.5\,\%$ lower despite the $22\,\%$ increase in
$B^{3}$. A fine entry places less waterplane area at large distance from the centerline, and near the upright condition this distribution governs the initial stability. The same fineness that reduces the wetted surface therefore also reduces the second moment of the waterplane. The reward of the generative framework contains a friction term proportional to the wetted surface but no term for the waterplane inertia, so that a reduction of wetted surface obtained in this way would be credited by the reward while the accompanying loss of initial stability would not be registered.

This example shows that a feasibility screen and an objective term are not equivalent means of enforcing a design requirement. A screen establishes an admissible region, and an optimization that has no gradient with respect to stability may reduce the stability margin up to the limit of that region without the reduction being reflected in the objective value. The large-angle
behavior of the two hulls further shows that their ranking depends on which stability measure is considered, since GenHull~E is less stiff near upright yet attains a larger maximum righting lever at a greater angle of heel. The
choice of which requirements are enforced as constraints, which are included in the objective, and which are screened after generation therefore has consequences that the objective value alone does not reveal.

\subsection{Requirements for a CFD-Based Reward}
\label{subsec:reward_design}
The result of Section~\ref{subsec:representation_limits} might suggest that improved training data would enable the generative model to recover a satisfactory bulbous bow. However, two distinct limitations are involved: geometric expressivity and reward fidelity.
The first is representational. The adopted parameterization cannot express a faired bulbous bow; consequently, no amount of training data can recover geometry outside the admissible design space. Addressing this limitation requires a richer parameterization or geometric representation.
The second concerns the physics represented by the reward. The present reward combines ITTC-57 friction based on wetted surface with an inviscid wave-resistance estimate and therefore omits viscous-pressure effects. These effects contribute to the penalty of the parameterized bulb (Section~\ref{subsec:representation_limits}) and the separated flow around GenHull~E (Section~\ref{subsubsec:genhull_limitations}). Addressing this limitation requires a higher-fidelity reward, for example one trained on viscous CFD data.
These limitations are coupled. A richer parameterization with the present reward could generate new geometries but rank them incorrectly, whereas a higher-fidelity reward with the present parameterization cannot recover geometries outside its admissible family. Reliable generative optimization therefore requires geometric expressivity and reward fidelity to advance together.

Two further properties of the reward bound the confidence that can be placed in the empirical ranking. The frictional term uses the wetted surface computed by the hull-generation routines of \cite{bagazinski2023ship}, which integrate the waterline perimeter over height and therefore undercount any part of the
hull surface that is not vertical. For the KCS-like parameter set from which GenHull~E was derived, this routine returns $45\,\%$ less wetted surface at the design draft than the area of the hull surface. The weight of this error in the reward depends on the share of friction in the aggregated resistance, which is small at the highest speeds that dominate the aggregate but large within the speed range of present interest. In addition, the reward contains no displacement term, and displacement is accounted for only in the matched-displacement screening of Section~\ref{subsec:empirical_screening}.

A resistance dataset generated with the pipeline of
Section~\ref{subsec:cfd_pipeline} would provide three properties that the
present reward lacks. Such a dataset contains viscous pressure, so that
bulb-like features can be scored; it is resolved by speed using a consistent
method, so that a crossover of the kind reported in
Section~\ref{subsubsec:genhull_resistance} is represented and not averaged
away; and, being produced under a single rule set, it permits hull-to-hull
differences to be interpreted to the extent established in
Section~\ref{subsubsec:cfd_validation}. Absolute accuracy is not among these
properties. The comparison envelope extends from $+14\,\%$ to $-9\,\%$ across
the full speed range and is approximately $\pm7\,\%$ within the range of
present interest, and a surrogate trained on such data would inherit this
bias. The acceptability of this bias rests on its being common-mode within a
hull family, which the two-hull validation supports and bounds at four to five
percentage points across dissimilar forms. The generation of such a dataset at
the scale required by a surrogate, of the order of 200 hulls at six speeds, is
the motivation for automating the evaluation pipeline and is the subject of
continuing work.

%% file: draft_v1.bbl
\begin{thebibliography}{10}
\expandafter\ifx\csname url\endcsname\relax
  \def\url#1{\texttt{#1}}\fi
\expandafter\ifx\csname urlprefix\endcsname\relax\def\urlprefix{URL }\fi
\expandafter\ifx\csname href\endcsname\relax
  \def\href#1#2{#2} \def\path#1{#1}\fi

\bibitem{gaspar_2012}
H.~M. Gaspar, D.~H. Rhodes, A.~M. Ross, S.~O. Erikstad, \href{https://doi.org/10.5957/jspd.2012.28.4.145}{Addressing complexity aspects in conceptual ship design: A systems engineering approach}, Journal of Ship Production and Design 28~(4) (2012) 145--159.

\bibitem{zhang2019research}
B.-J. Zhang, S.-L. Zhang, \href{https://doi.org/10.1007/978-981-10-8423-2}{Research on Ship Design and Optimization Based on Simulation-Based Design ({SBD}) Technique}, Springer, Singapore, 2019.

\bibitem{perez2011constrained}
F.~P{\'e}rez, J.~Clemente, Constrained design of simple ship hulls with b-spline surfaces, Computer-Aided Design 43~(12) (2011) 1829--1840.

\bibitem{han2012hydrodynamic}
S.~Han, Y.-S. Lee, Y.~B. Choi, \href{https://doi.org/10.1007/s00773-011-0148-8}{Hydrodynamic hull form optimization using parametric models}, Journal of Marine Science and Technology 17 (2012) 1--17.

\bibitem{lu2016hydrodynamic}
Y.~Lu, X.~Chang, A.-k. Hu, \href{https://doi.org/10.1080/19942060.2016.1159987}{A hydrodynamic optimization design methodology for a ship bulbous bow under multiple operating conditions}, Engineering Applications of Computational Fluid Mechanics 10~(1) (2016) 330--345.

\bibitem{holtrop1982approximate}
J.~Holtrop, G.~G.~J. Mennen, \href{https://doi.org/10.3233/ISP-1982-2933501}{An approximate power prediction method}, International Shipbuilding Progress 29~(335) (1982) 166--170.

\bibitem{coppede2019hydrodynamic}
A.~Copped{\`e}, S.~Gaggero, G.~Vernengo, D.~Villa, \href{https://doi.org/10.1016/j.apor.2019.05.026}{Hydrodynamic shape optimization by high fidelity {CFD} solver and {Gaussian} process based response surface method}, Applied Ocean Research 90 (2019) 101841.

\bibitem{harries2019faster}
S.~Harries, C.~Abt, \href{https://doi.org/10.1016/j.oceaneng.2019.106470}{Faster turn-around times for the design and optimization of functional surfaces}, Ocean Engineering 193 (2019) 106470.

\bibitem{khan2023shiphullgan}
S.~Khan, K.~Goucher-Lambert, K.~Kostas, P.~Kaklis, \href{https://doi.org/10.1016/j.cma.2023.116051}{{ShipHullGAN}: A generic parametric modeller for ship hull design using deep convolutional generative model}, Computer Methods in Applied Mechanics and Engineering 411 (2023) 116051.

\bibitem{sohl2015deep}
J.~Sohl-Dickstein, E.~Weiss, N.~Maheswaranathan, S.~Ganguli, \href{https://proceedings.mlr.press/v37/sohl-dickstein15.html}{Deep unsupervised learning using nonequilibrium thermodynamics}, in: Proceedings of the 32nd International Conference on Machine Learning, Vol.~37 of Proceedings of Machine Learning Research, PMLR, 2015, pp. 2256--2265.

\bibitem{ho2020denoising}
J.~Ho, A.~Jain, P.~Abbeel, \href{https://arxiv.org/abs/2006.11239}{Denoising diffusion probabilistic models}, in: Advances in Neural Information Processing Systems, Vol.~33, 2020, pp. 6840--6851.

\bibitem{bagazinski2023shipgen}
N.~J. Bagazinski, F.~Ahmed, \href{https://doi.org/10.3390/jmse11122215}{{ShipGen}: A diffusion model for parametric ship hull generation with multiple objectives and constraints}, Journal of Marine Science and Engineering 11~(12) (2023) 2215.

\bibitem{bagazinski2024c}
N.~J. Bagazinski, F.~Ahmed, \href{https://doi.org/10.59490/imdc.2024.841}{{C-ShipGen}: a conditional guided diffusion model for parametric ship hull design}, in: International Marine Design Conference, 2024.

\bibitem{krishnamoorthy2023diffusion}
S.~Krishnamoorthy, S.~M. Mashkaria, A.~Grover, \href{https://proceedings.mlr.press/v202/krishnamoorthy23a.html}{Diffusion models for black-box optimization}, in: Proceedings of the 40th International Conference on Machine Learning, Vol. 202 of Proceedings of Machine Learning Research, PMLR, 2023, pp. 17842--17857.

\bibitem{keramati2025heatgen}
H.~Keramati, M.~Sadeghi, R.~K. Jaiman, \href{https://arxiv.org/abs/2511.09578}{{HeatGen}: A guided diffusion framework for multiphysics heat sink design optimization}, arXiv preprint arXiv:2511.09578 (2025).

\bibitem{keramati2026reward}
H.~Keramati, P.~Kirchen, M.~Hannan, R.~K. Jaiman, \href{https://doi.org/10.1016/j.engappai.2025.113378}{A reward-directed diffusion framework for generative design}, Engineering Applications of Artificial Intelligence 165 (2026) 113378.

\bibitem{li2024derivative}
X.~Li, Y.~Zhao, C.~Wang, G.~Scalia, G.~Eraslan, S.~Nair, T.~Biancalani, S.~Ji, A.~Regev, S.~Levine, M.~Uehara, \href{https://arxiv.org/abs/2408.08252}{Derivative-free guidance in continuous and discrete diffusion models with soft value-based decoding}, arXiv preprint arXiv:2408.08252 (2024).

\bibitem{del2006sequential}
P.~Del~Moral, A.~Doucet, A.~Jasra, \href{https://doi.org/10.1111/j.1467-9868.2006.00553.x}{Sequential {Monte Carlo} samplers}, Journal of the Royal Statistical Society Series B: Statistical Methodology 68~(3) (2006) 411--436.

\bibitem{shields2016generalization}
M.~D. Shields, J.~Zhang, The generalization of latin hypercube sampling, Reliability Engineering \& System Safety 148 (2016) 96--108.

\bibitem{wandb}
L.~Biewald, \href{https://www.wandb.com/}{Experiment tracking with weights and biases}, software available from wandb.com (2020).

\bibitem{hirt1981vof}
C.~W. Hirt, B.~D. Nichols, Volume of fluid ({VOF}) method for the dynamics of free boundaries, Journal of Computational Physics 39~(1) (1981) 201--225.

\bibitem{menter2003sst}
F.~R. Menter, M.~Kuntz, R.~Langtry, Ten years of industrial experience with the {SST} turbulence model, Turbulence, Heat and Mass Transfer 4~(1) (2003) 625--632.

\bibitem{tokyo2015kcs}
{Tokyo 2015 Workshop}, \href{https://t2015.nmri.go.jp/Instructions_KCS/instruction_KCS.html}{{KRISO} container ship ({KCS}): {EFD} data and format of computed results, case 2.1}, Tokyo 2015: A Workshop on {CFD} in Ship Hydrodynamics, National Maritime Research Institute (2015).

\bibitem{kim2001kcs}
W.~J. Kim, S.~H. Van, D.~H. Kim, Measurement of flows around modern commercial ship models, Experiments in Fluids 31~(5) (2001) 567--578.

\bibitem{elmoctar2012dtc}
O.~el~Moctar, V.~Shigunov, T.~Zorn, {Duisburg Test Case}: Post-panamax container ship for benchmarking, Ship Technology Research 59~(3) (2012) 50--64.

\bibitem{feng2021parametric}
Y.~Feng, O.~el~Moctar, T.~E. Schellin, Parametric hull form optimization of containerships for minimum resistance in calm water and in waves, Journal of Marine Science and Application 20~(4) (2021) 670--693.

\bibitem{ittc2021verification}
{International Towing Tank Conference}, Uncertainty analysis in {CFD} verification and validation methodology and procedures, Tech. Rep. 7.5-03-01-01, ITTC (2021).

\bibitem{keramati2026parametric}
H.~Keramati, M.~A. Hannan, P.~Kirchen, S.~Tousi, J.~Riekkinen, M.~Barrett, J.~Pyra, R.~K. Jaiman, A parametric diffusion framework for localized generative optimization of ship hulls, in: International Conference on Offshore Mechanics and Arctic Engineering, Vol. 89541, American Society of Mechanical Engineers, 2026, p. V006T08A050.

\bibitem{bagazinski2023ship}
N.~J. Bagazinski, F.~Ahmed, {Ship-D}: Ship hull dataset for design optimization using machine learning, in: International Design Engineering Technical Conferences and Computers and Information in Engineering Conference, Vol. 87301, American Society of Mechanical Engineers, 2023, p. V03AT03A028.

\end{thebibliography}
